\magnification=\magstep1
\input amstex
\documentstyle{amsppt}
\catcode`\@=11 \loadmathfont{rsfs}
\def\mycal{\mathfont@\rsfs}
\csname rsfs \endcsname \catcode`\@=\active

\vsize=6.5in

\topmatter 
\title 
Smooth bimodules \\ and cohomology of II$_1$ factors 
\endtitle
\author  Alin Galatan and Sorin Popa  \endauthor

\rightheadtext{Smooth bimodules and cohomology }

\affil     {\it  University of California, Los Angeles} \endaffil

\address Math.Dept., UCLA, Los Angeles, CA 90095-1555\endaddress
\email  agalatan\@math.ucla.edu, popa\@math.ucla.edu\endemail

\thanks S.P. supported in part by NSF Grant DMS-1400208 \endthanks

\abstract  We prove that, under rather general conditions, the 1-cohomology 
of a von Neumann algebra $M$ with values in a Banach $M$-bimodule satisfying a combination of smoothness and 
operatorial conditions, vanishes. For instance, we show that if $M$ acts normally on a Hilbert space $\Cal H$ and 
$\Cal B_0\subset \Cal B(\Cal H)$ is a norm closed $M$-bimodule such that any $T\in \Cal B_0$ is {\it smooth} (i.e. the left and right multiplication 
of $T$ by $x\in M$ are continuous from the unit ball of $M$ with the $s^*$-topology to $\Cal B_0$ with its norm), then any derivation of $M$ into $\Cal B_0$ is inner. 
The compact operators are smooth over any $M\subset \Cal B(\Cal H)$, but there is a large variety of non-compact smooth elements as well.

\endabstract

\endtopmatter

\document

\heading 0. Introduction \endheading

Given a von Neumann algebra $M$ and a Banach $M$-bimodule $\Cal B$, an element $T\in \Cal B$ is {\it smooth} (over $M$) if 
the left and right multiplication of $T$ by elements in $M$ are continuous from the unit ball of $M$ 
with the $s^*$-topology, into $\Cal B$ with its norm topology. The space $s_M^*(\Cal B)$ (or $s^*(\Cal B)$  for simplicity) of all smooth elements in $\Cal B$ is itself a 
Banach $M$-bimodule and 
we investigate in this paper the  question of whether the 1st cohomology with values in a smooth closed $M$-bimodule 
$\Cal B_0\subset s^*(\Cal B)$ vanishes. In other words, whether any derivation $\delta: M \rightarrow \Cal B_0$ 
(i.e., a linear map satisfying $\delta(xy)=x\delta(y)+\delta(x)y$, $\forall x,y\in M$) 
is ``inner'', in the sense that there exists $T\in \Cal B_0$ such that $\delta(x)=Tx-xT$, $\forall x\in M$. 

We in fact only consider Banach $M$-bimodules $\Cal B$ that 
are {\it operatorial} (over $M$), i.e., for which the norm on $\Cal B$ satisfies the axiom  $\|pTp+(1-p)T(1-p)\|=\max \{\|pTp\|, 
(1-p)T(1-p)\|\}$, for all $T\in \Cal B$ and all projections $p$ in $M$, and they will usually be   
assumed dual and normal. 

The prototype example of a dual normal operatorial $M$-bimodule is the space of all linear bounded operators $\Cal B(\Cal H)$ 
on the Hilbert space $\Cal H$ on which 
$M$ is represented.  Its smooth part $s^*(\Cal B(\Cal H))$ is a hereditary $C^*$-subalgebra of $\Cal B(\Cal H)$ that contains 
the space of compact operators $\Cal K(\Cal H)$.  But  unless $M$ is a direct sum of matrix algebras, $s^*(\Cal B(\Cal H))$ 
 is much bigger, containing  a large variety of smooth elements that are non-compact (see 1.6, 1.7 hereafter).  

We first prove that if $\Cal B$ is a dual normal operatorial $M$-bimodule, 
then any derivation $\delta$ from  $M$ into a closed sub-bimodule $\Cal B_0\subset s^*(\Cal B)$ 
can be ``integrated'' to an element in $\Cal B_0$ on any abelian von Neumann subalgebra of $M$,  
and show that if $M$ satisfies some very weak regularity condition (e.g. if $M$ has a Cartan subalgebra, or if it has the property $(\Gamma)$ of [MvN2]), then 
any element in $\Cal B_0$ implementing $\delta$ on a diffuse abelian subalgebra 
automatically implements $\delta$ on all $M$. We also prove a similar statement for 
{\it smooth derivations} of $M$ into an arbitrary (not necessarily dual normal) 
operatorial Banach bimodule $\Cal B_0$, i.e. for derivations that are continuous from the unit ball of $M$ with the $s^*$-topology into 
$\Cal B_0$ with its norm topology. The precise statement is as follows: 

\proclaim{0.1. Theorem} Let $M$ be a von Neumann algebra with the property that  any $\text{\rm II}_1$ factor summand of $M$ contains   
a wq-regular diffuse abelian von Neumann subalgebra. Let $\Cal B_0$ be an operatorial Banach $M$-bimodule and $\delta:M \rightarrow \Cal B_0$ a derivation. 
Assume that either $\delta$ is smooth, or $\Cal B_0$ is a closed sub-bimodule of the smooth part of a dual normal operatorial $M$-bimodule $\Cal B$, 
i.e. $\Cal B_0\subset s^*(\Cal B)$. Then  
there exists $T\in \Cal B_0$ such that $\delta=\text{\rm ad}T$,  $\|T\|\leq \|\delta\|$.  
\endproclaim

The {\it wq-regularity} condition in the above statement requires that $M$ has a diffuse abelian von Neumann subalgebra $A\subset M$ such that 
one can ``reach out'' from $A$ to $M$ inductively, by a chain of algebras, $A=N_0 \subset N_1 \subset .... \subset N_{\imath} =M$, 
such that for all $\jmath < \imath$, $N_{\jmath +1}$ is generated by unitaries $u\in M$ with 
the property that $uN_{\jmath}u^*\cap N_{\jmath}$ is diffuse. We mention that in the case a derivation 
$\delta$ takes values into a smooth $M$-bimodule $\Cal B_0$ that's contained in a dual normal $M$-bimodule,  
we actually prove that $\delta$ follows automatically smooth, without the wq-regularity condition. 

Our second main result shows that in the specific case $\Cal B=\Cal B(\Cal H)$, 0.1 holds true even without the wq-regularity assumption: any smooth valued 
derivation of any von Neumann subalgebra $M$ into $\Cal B(\Cal H)$ is inner. More precisely, we prove:

\proclaim{0.2. Theorem}  Let $M\subset \Cal B(\Cal H)$ be an arbitrary  von Neumann algebra represented on a Hilbert space $\Cal H$. Let $\Cal B_0\subset s^*(\Cal B(\Cal H))$ be a 
Banach $M$ sub-bimodule and $\delta: M \rightarrow \Cal B_0$ a derivation. Then there exists $T\in \Cal B_0$ such that $\delta=\text{\rm ad}T$, $\|T\|\leq \|\delta\|$. 
\endproclaim

In particular, since the compact operators are smooth over any von Neumann algebra, the above result recovers a result in  [P3], 
sowing that any derivation of a von Neumann algebra with values into the space of compact operators, is implemented by a compact operator. 
It also provides an answer to a question posed by Gilles Pisier, concerning derivations  from a C$^*$-algebra with a trace $(M_0, \tau)$ 
(such as the reduced C$^*$-algebra of a free group, $C^*_r(\Bbb F_n)$, $2\leq n \leq \infty$), represented on a Hilbert space $\Cal H$, 
with values in $\Cal B(\Cal H)$, satisfying $\| \ \|_2$ - $\| \ \|_{\Cal B(\Cal H)}$ type continuity conditions. Indeed, 
the existence of such a derivation $\delta$ implies that the  representation of $M_0$ is automatically ``normal with respect to $\tau$'' 
on the non-degenerate part of $\delta$, 
and that $\delta$ extends to a smooth derivation of $M=M_0''$ into $s^*_M(\Cal B(\Cal H))$, so 0.2 applies to get: 

\proclaim{0.3. Corollary}  Let $M_0$ be a $C^*$-algebra with a faithful trace $\tau$ and 
$M_0\subset \Cal B(\Cal H)$ a faithful representation of $M_0$. Let $\delta:M_0 \rightarrow \Cal B(\Cal H)$ be a derivation. Assume $\delta$ is continuous from the unit 
ball of $M_0$ with the topology given by the Hilbert norm $\|x\|_2=\tau(x^*x)^{1/2}$, $x\in M_0$, 
to $\Cal B(\Cal H)$ with the operator norm topology. Then  
there exists $T\in \Cal B(\Cal H)$ such that $\delta=\text{\rm ad}T$ and $\|T\|\leq \|\delta\|$.  \endproclaim

As we mentioned before, the prototype example of a smooth operatorial $M$-bimodule is the space of compact operators $\Cal K(\Cal H)$ on the Hilbert space $\Cal H$ on which 
$M$ is represented, a case studied in [JP], [P2,3]. The initial motivation for our work  has in fact been to provide an abstract setting 
for this case, and find the largest degree of generality for which arguments in the spirit of  [JP], [P2,3] can be carried over. At the same time,  
we were hoping to find $M$-bimodules for which the vanishing cohomology on abelian subalgebras (as in [JP]) and automatic extension 
properties (as in [P2]) do hold,  while the  arguments in [P3], showing the vanishing cohomology for arbitrary II$_1$ factors, don't. 

Results that show vanishing of the 1-cohomology 
with coefficients in $\Cal B(\Cal H)$ for arbitrary II$_1$ factors, like in Theorem 0.2 above, can be relevant for the similarity problem (cf. [Ch], [Ki], [Pi3]). In turn,  
the existence of a smooth $M$-bimodule $\Cal B_0$ for which a II$_1$ factor $M$  
has non-inner derivations, would imply, via Theorem 0.1, that $M$ has no 
diffuse abelian wq-regular subalgebra (so  in particular, $M$ would have no Cartan subalgebras, would be prime, etc). 
This falls within the larger scope of finding a cohomology theory for II$_1$ factors that's non-vanishing (and if possible calculable) and that could detect important properties 
of II$_1$ factors, such as 
absence of regularity, or infinite generation. We do not provide smooth such $\Cal B_0$'s here and it may be that in fact 
0.1, 0.2 hold true for all II$_1$ factors and all smooth bimodules. We leave the clarification of this aspect as an open problem.

The paper is organized as follows. In Section 1 we define the smooth part of  a Banach $M$-bimodule, 
prove the main properties (1.2, 1.3, 1.8) and provide examples. In particular, we show that if $M\subset \Cal B(\Cal H)$ is a diffuse finite von Neumann algebra  then  
the smooth part of the $M$-bimodule $\Cal B(\Cal H)$ contains infinite dimensional projections along ``Gaussian'', ``CAR'' and ``free'' directions (see 1.6, 1.7). 
In Section 2 we consider smooth valued derivations and   
show that they can be ``integrated'' on abelian subalgebras (see 2.5). The proof of this result follows ideas and  techniques from  [JP]. 
In Section 3 we introduce the notion 
of wq-regularity for subalgebras and prove Theorem 0.1 (see 3.7, 3.8), 
by showing that if a smooth valued derivation $\delta$ has been integrated to a smooth element $T$ on a diffuse abelian subalgebra that's wq-regular in $M$, then $T$ automatically 
implements  $\delta$ on all $M$. We in fact prove a stronger result, involving a property that generalizes properties $(\Gamma)$ in [MvN2] and $(C)$ in [P2], 
and that we call $(C')$ (see 3.6).   
The proofs in this section consist of a refinement of 
arguments in [P2]. 

In Section 4 we prove Theorem 0.2 and deduce Corollary 0.3 (see 4.1 and 4.5). The proof of 4.1 follows the same strategy as in the case of compact valued derivations 
in [P3],  but the additional technical difficulties are significant. We overcome them 
by using the incremental patching technique in [P4,6,7].  Section 5 contains general remarks, 
including a generalization of Theorem 4.1 that recovers results in  [PR] (see 5.1),  the definition of smooth $n$-cohomology for $n\geq 2$ (5.2) and 
some final comments on smooth cohomology (5.3). 
For convenience, we have included an Appendix with the proof of a general ``continuity principle'', 
extracted from [P2,3], used several times in this paper. For the basics of von Neumann algebras, we refer the reader to the classic monographs [D], [T].

\heading 1. Smooth bimodules \endheading

Recall that if $M$ is a unital Banach algebra (which we will in fact always assume in this paper 
to be a von Neumann algebra), then a {\it Banach $M$-bimodule} $\Cal B$ is an $M$-bimodule with 
the property that $1_MT=T1_M=T$, $\forall T\in \Cal B$, and the left and right multiplication operations 
$M \times \Cal B \ni (x, T) \mapsto xT\in  \Cal B$, $\Cal B \times M \ni (T, x) \mapsto Tx \in \Cal B$ are bounded bilinear maps. 
For simplicity, we will in fact always assume that $\|xT\|\leq \|x\| \|T\|$, $\|Tx\|\leq \|T\| \|x\|$, $\forall x\in M, T\in \Cal B$. 

If in addition $\Cal B$ is the dual of a Banach space $\Cal B_*$ and for each $x\in M$ the maps $\Cal B \ni T \mapsto xT \in \Cal B$, 
$\Cal B \ni T \mapsto Tx \in \Cal B$ are continuous with respect to the $\sigma(\Cal B, \Cal B_*)$ topology 
(also called weak$^*$-topology), then 
$\Cal B$ is called a {\it dual  $M$-bimodule}.  Finally, if $M$ is a von Neumann algebra, $\Cal B$ is a dual $M$-bimodule and for each $T\in \Cal B$ 
the maps $M \ni x \mapsto xT \in \Cal B$, $M\ni T \mapsto Tx \in \Cal B$  are continuous from $(M)_1$ with the $\sigma(M, M_*)$-topology 
(also called ultraweak topology), then we say that the 
dual $M$-bimodule $\Cal B$ is {\it normal}. 
\vskip .1in 

\noindent
{\bf 1.1. Definitions.} $1^\circ$ Let $M$ be a von Neumann algebra and $\Cal B$ a Banach $M$-bimodule. An element $T\in \Cal B$ 
is {\it smooth} (with respect to $M$) if the maps $x \mapsto xT$, $x \mapsto T x$ are continuous 
from $(M)_1$ with the $s^*$-topology to $\Cal B$ with the norm topology. 
We denote by 
$s_M^*(\Cal B)$ (or simply $s^*(\Cal B)$ if no confusion is possible), 
the set of smooth elements  $T \in \Cal B$ and call it the {\it smooth part} of the $M$-bimodule $\Cal B$.  
The Banach $M$-bimodule $\Cal B$ 
is {\it smooth} if $s^*(\Cal B)=\Cal B$.  A subset $\Cal S$ of a Banach $M$-bimodule $\Cal B$ is {\it uniformly smooth} if $\forall \varepsilon > 0$, 
there exists a $s^*$-neighborhood $\Cal V$ of $0$ in $M$ such that if $x\in (M)_1 \cap \Cal V$, then $\|T x\|, \|x T \|\leq \varepsilon$, $\forall T \in \Cal S$. 
Note that if $(M,\tau)$ is a finite von Neumann algebra with a faithful normal trace $\tau$ (the case of most interest for us), then this amounts to the existence of some 
$\alpha > 0$ such that if $x\in (M)_1$, $\|x\|_2 \leq \alpha$, then $\|xT\|, \|Tx\|\leq \varepsilon$, where as usual $\|x\|_2=\tau(x^*x)^{1/2}$ denotes the Hilbert norm 
implemented by $\tau$.

$2^\circ$  An $M$-bimodule is {\it operatorial} if for any $p\in \Cal P(M)=\{p\in M \mid p=p^*=p^2\}$  and any $T \in \Cal B$, we have $\|pT p + (1-p)T (1-p)\|
=\max \{\|pT p\|, \|(1-p)T (1-p)\|\}$.

\proclaim{1.2. Proposition} Let $\Cal B$ be a Banach $M$-bimodule and $\Cal B_0\subset \Cal B$ a Banach sub-bimodule. 

$1^\circ$  $s^*(\Cal B)$ is a Banach $M$-bimodule. 

$2^\circ$ If $M_0\subset M$ is a von Neumann subalgebra, then $s^*_M(\Cal B)\subset s^*_{M_0}(\Cal B)$. 

$3^\circ$ $s^*(\Cal B)/\Cal B_0\subset s^*(\Cal B/\Cal B_0)$.

$4^\circ$ If $\Cal B$ is operatorial, then $\Cal B/\Cal B_0$, with its quotient norm and $M$-bimodule structure, is an operatorial 
$Banach$ $M$-bimodule. 

$5^\circ$ The bidual $\Cal B^{**}$ of $\Cal B$ has a natural dual $M$-bimodule structure and if $\Cal B$ is operatorial, then so is $\Cal B^{**}$. 

$6^\circ$ Let $\Cal B$ be a dual normal $M$-bimodule. Then its predual $\Cal B_*$ has a natural Banach $M$-bimodule structure which is smooth $($ i.e.  
$s^*_M(\Cal B_*)=\Cal B_*)$. But if $\Cal B$ is operatorial, then the norm on $\Cal B_*$ satisfies $\|\varphi(p  \cdot p) + \varphi((1-p) \cdot (1-p))\|=
\|\varphi(p  \cdot p)\| + \|\varphi((1-p) \cdot (1-p))\|$, for $\varphi\in \Cal B_*$, $p\in \Cal P(M)$, and thus it is not operatorial in general.  
\endproclaim 
\noindent
{\it Proof}. $1^\circ$ $s^*(\Cal B)$ is clearly closed to addition and scalar multiplication. Also, if $T_n \in s^*(\Cal B)$ and $\lim_n \|T_n-T\|=0$ 
for some $T\in \Cal B$, then for any $\varepsilon > 0$ there exists $n_0$ such that $\|T-T_{n_0}\|\leq \varepsilon/2$. 
Also, since $T_{n_0}\in s^*(\Cal B)$, there exists an $s^*$-neighborhood $\Cal V$ of $0$ in $M$ 
such that if $x\in M \cap \Cal V$ then $\|xT_{n_0}\|, \|T_{n_0}x\|\leq \varepsilon/2$. Thus, for such $x$ we also have 
 $\|xT\| \leq \|xT_{n_0}\| + \|x(T-T_{n_0})\| \leq \varepsilon$. This shows that $T\in s^*(\Cal B)$. 
 
If $\|xT\|, \|Tx\|\leq \varepsilon$ for $x\in \Cal V \cap (M)_1$, for some $s^*$-neighborhood $\Cal V$ of $0$ in $(M)_1$, 
then in particular  it holds true for $x\in \Cal V \cap (M_0)_1$, for any subalgebra $M_0\subset M$. Moreover, $\|xT\|_{\Cal B/\Cal B_0}, \|Tx\|_{\Cal B/\Cal B_0}\leq \varepsilon$ 
as well, proving $2^\circ$ and $3^\circ$. 

Parts $4^\circ, 5^\circ$ are trivial by the definitions. The fact that $s^*_M(\Cal B_*)=\Cal B_*$ in $6^\circ$ is (Lemma 5 in [R]) and the last part of $6^\circ$ is trivial.  

\hfill $\square$

\proclaim{1.3. Proposition} $1^\circ$ If $M$ is embedded as a $^*$-subalgebra in a unital $C^*$-algebra $\Cal B$ with $1_M=1_{\Cal B}$, 
then $\Cal B$ with its left and right multiplication by elements in $M$ is an operatorial Banach $M$-bimodule 
and $s_M^*(\Cal B)$ is a hereditary $C^*$subalgebra of $\Cal B$ which is both an $M$-bimodule and an $M'\cap \Cal B$-bimodule. 
\vskip .05in 

$2^\circ$ Assume  $M\subset \Cal B(\Cal H)$ is a normal representation of $M$. Then the hereditary $C^*$-algebra 
$s^*(\Cal B(\Cal H))\subset \Cal B(\Cal H)$  
contains the space of compact operators $\Cal K(\Cal H)$ and it is both an $M$-bimodule and an $M'$-bimodule. 

$3^\circ$ With $M\subset \Cal B(\Cal H)$ as in $2^\circ$ let $P\subset M$ be a von Neumann subalgebra 
and $\Cal H_0\subset \Cal H$ a Hilbert subspace such that the projection $p$ of $\Cal H$ onto $\Cal H_0$ 
commutes with $P$. Then $ps^*_M(\Cal B(\Cal H))p \subset  s_{Pp}^*(\Cal B(\Cal H_0))$, with equality if $P=M$. 

$4^\circ$ If $(M,\tau)$ is a finite von Neumann algebra, $P\subset M$ is a von Neumann subalgebra and $e_P$ 
denotes the orthogonal projection of $L^2M$ onto $L^2P$, then for any $T\in s_M^*(\Cal B(L^2M))$ we have $e_PT, Te_p 
\in s^*_M(\Cal B(L^2M))$. Also,  with the usual identifications, we have $e_Ps^*_M(\Cal B(L^2M))e_P=s^*_P(\Cal B(L^2P))$. 
\endproclaim 
\noindent
{\it Proof}. $1^\circ$ If $T_1, T_2 \in s^*(\Cal B)$ and $x_i\in (M)_1$ converges  $s^*$ to $0$, 
then $\|x_n(T_1T_2)\| \leq \|x_iT_1\| \|T_2\| \rightarrow 0$ and $\|(T_1T_2)x_i\|\leq \|T_1\| \|T_2 x_i\|\rightarrow 0$, 
showing that $T_1T_2 \in s^*(\Cal B)$. Also, $\|x_i T_1^*\|=\|T_1x_i^*\| \rightarrow 0$, thus $s^*(\Cal B)$ is actually a $C^*$-subalgebra 
of $\Cal B$. Moreover, if $T_0 \in \Cal B$ with $0\leq T_0\leq T\in s^*(\Cal B)$, then $\|x_iT_0\|\leq \|x_iT\|\rightarrow 0$ and 
so $s^*(\Cal B)$ is hereditary. 

The smoothness of an element is clearly  invariant to 
left-right multiplication by elements commuting with $M$.

$2^\circ$ We in fact have that for any bounded net of operators $x_i\in \Cal B(\Cal H)$ converging $s^*$ to $0$, 
and any $T\in \Cal K(\Cal H)$, $\lim_i \|x_iT\| = \lim_i \|Tx_i\|=0$, so in particular $\Cal K(\Cal H)\subset s^*(\Cal B(\Cal H))$.  

$3^\circ$  follows trivially from the fact that smoothness over $P$ is invariant to left-right multiplication by elements in $P'$. 

$4^\circ$ If $x_n\in (M)_1$ satisfies $\tau(x_n^*x_n)=\|x_n\|_2^2 \rightarrow 0$, then by using that $e_P(y)=E_P(y)e_P$, 
where $E_P$ denotes the $\tau$-preserving conditonal expectation of $M$ onto $P$, and the fact 
that $E_P(x_n^*x_n)$ converges to $0$ in the $s^*$-topology, we get: 
$$
\|x_ne_PT\|^2=
\|T^*e_Px_n^*x_ne_PT\|=\|T^*e_PE_P(x_n^*x_n)T\| \leq \|T\| \|E_P(x_n^*x_n)T\|\rightarrow 0.  
$$
Similarly $\|Te_Px_n\|\rightarrow 0$.  The same calculation shows that if $T_0\in s^*_P(\Cal B(L^2P))$, then 
$e_PT_0e_P$ viewed as an element in $\Cal B(L^2M)$ is smooth$/M$. Combined with $3^\circ$, this shows that 
$e_Ps^*_M(\Cal B(L^2M))e_P=s^*_P(\Cal B(L^2P))$. 
\hfill 
$\square$

\proclaim{1.4. Lemma} Let $(M,\tau)$ be a finite von Neumann algebra,  $M\subset \Cal B(\Cal H)$ a normal representation of $M$,  
$\Cal H_0\subset \Cal H$ a Hilbert subspace and $p=\text{\rm proj}_{\Cal H_0}$.   

$1^\circ$ $p \in  s^*(\Cal B(\Cal H))$ iff   $\lim_{\alpha \rightarrow 0} 
\sup \{\| x(\xi) \| \mid x \in (M)_1, \|x\|_2 \leq \alpha, \xi \in (\Cal H_0)_1\} =0$ and iff $\lim_n \|q_npq_n\|=0$  $($equivalently, $\|(1-q_n)p(1-q_n) - p\|\rightarrow 0)$, 
for any sequence of projections $\{q_n\}_n \subset \Cal P(M)$ that decrease to $0$.  

$2^\circ$ If $\Cal H=L^2M$, then a sufficient condition for $p$ to be in $s^*(\Cal B(L^2M))$  
is that $\Cal H_0 \subset \hat{M}$ and $\exists C>0$ such that $\|x\| \leq C \|x\|_2$, $\forall x\in \Cal H_0$. 
\endproclaim 
\noindent
{\it Proof}. The first equivalence in $1^\circ$ is  just a reformulation of the condition of smoothness for the orthogonal projection $\text{\rm proj}_{\Cal H_0}$. 
If $q_n$ are mutually orthogonal projections in $\Cal P(M)$, they converge to $0$ in the strong operator topology, so $\|q_np\|\rightarrow 0$, as 
required. Conversely, if this latter condition is satisfied for all sequence of orthogonal projections, then by (proof of Theorem in [P1]; see Lemma A.1 here enclosed), 
it follows that $\|x_np\| \rightarrow 0$ for any sequence $\{x_n\}_n\subset (M)_1$ with $\tau(x_n^*x_n)\rightarrow 0$. 

$2^\circ$ If $\hat{x} \in \Cal H_0\subset \hat{M}\subset L^2M$, satisfies $\|x\|_2\leq 1$ and we take $y\in M$, then 
$\|y (\hat{x})\|_2=\|yx\|_2\leq \|y\|_2 \|x\| \leq C \|y\|_2$. This shows in particular that if $y\in (M)_1$ satisfies $\|y\|_2 \leq \alpha$ 
then $\|y p\| \leq C\alpha$, so $p$ is smooth. 
\hfill
$\square$

\vskip .05in
\noindent
{\bf 1.5. Definition} Let $(M,\tau)$ be a finite von Neumann algebra. An orthogonal projection $p$ 
of $L^2(M)$ onto a Hilbert subspace $\Cal H_0\subset L^2M$ is {\it strongly smooth} (or {\it s-smooth}) if the condition 1.4.2$^\circ$ is satisfied, 
i.e. there exists some $C>0$ such that $\|x\| \leq C\|x\|_2$, $\forall x\in \Cal H_0$. 

\vskip .05in

The next two results provide a large variety of concrete examples of smooth and strongly smooth non-compact elements 
in the case $M$ is a finite von Neumann algebra and the target bimodule is $\Cal B(L^2M)$. We distinguish three remarkable classes 
of infinite rank smooth projections $p=\text{\rm proj}_{\Cal H_0}$: 

\vskip .05in 

$(a)$ any ``Gaussian Hilbert subspace'' $\Cal H_0\subset L^2M$, along an arbitrary diffuse abelian von  
Neumann subalgebra $A$ of $M$ (see 1.6.2$^\circ$ below); 

$(b)$ any ``CAR Hilbert subspace'' $\Cal H_0\subset L^2M,$ along an arbitrary hyperfinite II$_1$ subfactor $R$ of $M$ (see $1.6.4^\circ$ 
below); 

$(c)$ any ``free Hilbert subspace'' $\Cal H_0\subset L^2M$, generated by unitaries in $\hat{M}\subset L^2M$ that are free independent (see 1.7.1$^\circ$).

\vskip .05in

All these spaces are  ``thin'' with respect to $M$, with the  abelian ``Gaussian directions'' being particularly thin (they are smooth but not strongly smooth). 
It is quite interesting that s-smoothness is a purely non-commutative phenomenon, that only occurs when $M$ is II$_1$.  

We are grateful to Assaf Naor and Gideon Schechtman for their help with the Gaussian example 1.6.2$^\circ$, and to Gilles Pisier for his help with the CAR example 1.6.4$^\circ$, below.

\proclaim{1.6. Proposition} Let $(M, \tau)$ be a finite von Neumann algebra and $M\subset \Cal B(L^2M)$ its standard representation.  

$1^\circ$ If $M$ is atomic, then $s^*(\Cal B(L^2M)) = \Cal K(L^2M)$. 

$2^\circ$ Assume $M$ is diffuse and let $A\subset M$ be a separable diffuse abelian von Neumann subalgebra. Let $\{\xi_n\}_n \subset L^2(A, \tau)$ 
be a sequence of independent and identically distributed real valued Gaussian random variables. If $p$ denotes the orthogonal projection of 
$L^2M$ onto the Hilbert subspace $\Cal H_0$ generated by $\{\xi_n\}_n$, then  $p\in s_M^*(\Cal B(L^2M))$. Thus, $\Cal K(L^2M) \subsetneq s^*(\Cal B(L^2M))$. 

$3^\circ$ If $M$ is diffuse abelian and $\Cal H=L^2M$, then any s-smooth projection $p\in \Cal B(L^2M)$ must have finite rank.   

$4^\circ$ Assume $M$ is type $\text{\rm II}_1$. If $\{u_n\}_n\subset \Cal U(M)$ is a sequence of selfadjoint unitaries of trace $0$ 
satisfying the canonical anticommutation relations $($CAR$)$ $u_iu_j=-u_ju_i$, $\forall i\neq j$, and we denote by $\Cal H_0 \subset L^2M$ 
the Hilbert space with orthonormal basis $\{u_n\}_n$, then the projection of $L^2M$ onto $\Cal H_0$ is s-smooth. 
\endproclaim 
\noindent
{\it Proof}. $1^\circ$ If $M$ has finite dimensional center, then $M$ is finite dimensional and the statement becomes trivial. 
If in turn $\Cal Z(M)$ is infinite dimensional, then it follows isomorphic to $\ell^\infty(\Bbb N)$ (because $M$ has a faithful trace). 
Let $\{e_k\}_{k\geq 1}$ be the atoms in $\Cal Z(M)$ and $p_n = \Sigma_{i=1}^ne_k$. Then $p_n$ are finite rank projections, i.e. $p_n\in \Cal K=\Cal K(L^2M)$,  
and they increase to $1$. Thus, if $T\in \Cal B=\Cal B(L^2M)$ is a smooth element, then $\|(1-p_n)T\|\rightarrow 0$, 
so in the Calkin algebra $\Cal B/\Cal K$ we have $\|T\|_{\Cal B/\Cal K}=\|(1-p_n)T\|_{\Cal B/\Cal K} \leq \|(1-p_n)T\| \rightarrow 0$, showing that $T\in \Cal K$. 

$2^\circ$ By $1.3.4^\circ$, it is sufficient to 
prove that if $p_n \in A$ are projections with $\tau(p_n)\rightarrow 0$, then $\sup \{ \|p_n\xi\|_2 \mid \xi=\xi^* \in (\Cal H_0)_1 \} \rightarrow 0$. 
Represent  $A$ as $L^\infty(\Bbb T)$ with $\tau$ corresponding to the  integral over the Haar measure  $\mu$. 
If $\xi = \Sigma_n c_n \xi_n$ with $c_n \in \Bbb R$, $\Sigma_n |c_n|^2=1$, 
then $\xi$ is still a Gaussian with the same distribution as the $\xi_n$'s. Since the $L^2$-norm of the restriction 
of a Gaussian $\xi$ to the set $X_n= \{t \mid |\xi(t)|\geq n\}$ decays exponentially in $n^2$, 
it follows that the $L^2$-norm of the restriction of $\xi$ to an arbitrary set $X\subset \Bbb T$ of measure $\mu(X) \leq \mu(X_n)$ is 
majorised by $\|\xi \chi_{X_n}\|_2$ and thus tends to $0$ as $n$ tends to $\infty$, independently of $\xi$. 

$3^\circ$ It is well known that if $\Cal H_0\subset L^\infty(\Bbb T)\subset L^2(\Bbb T)$ 
has dimension $1\leq n \leq \infty$ and satisfies the property that $\|x\|\leq C \|x\|_2$, $\forall x\in \Cal H_0$, 
then $C \geq \sqrt{n}$ (see e.g. [Pi1]). 

$4^\circ$ Since $M$ is type II$_1$,  it contains a copy of the hyperfinite II$_1$ factor $R$, which we can view as the weak closure of an 
increasing sequence of $2^n \times 2^n$ matrix algebras $B_n\subset R$. Moreover, by using the representation of 
$\cup_n B_n$ as the algebra of the canonical anticommutation relations (see e.g. [PoS]), we can find inside it  
selfadjoint unitaries $\{u_k\}_k$ with the properties: $(a)$ $u_1, ..., u_n \in B_n$ and they generate $B_n$ as an algebra; $(b)$ $\tau (u_k)=0$, $\forall k\geq 1$; 
$(c)$ $u_iu_j = -u_j u_i$, $\forall i\neq j$.   But then for any $c_k \in \Bbb R$ with $\Sigma_k |c_k|^2 < \infty$, 
the element $x=\Sigma_k c_ku_k$ satisfies  the identities:
$$
x^*x=(\Sigma_i c_i u_i)(\Sigma_j c_j u_j)=\Sigma_k c_k^2 + \Sigma_{i\neq j} c_ic_j u_i u_j
$$
$$
= \|x\|_2^2 1_M + \Sigma_{i<j} c_ic_j(u_iu_j + u_ju_i)= \|x\|_2^2 1_M. 
$$
Thus,  $x^*x = \|x\|_2^2 1_M$. This shows that if we denote by  $\Cal H_n$ 
the span of $\{u_k \mid 1\leq k \leq n\}$, then for any $x=x^*\in \Cal H_n$ we have $x^*x = \|x\|_2^2 1_M$, 
so in particular $\|x\| = \|x\|_2$. 
For arbitrary $x\in \Cal H_n$ we thus get  
$$
\|x\| \leq \|\Re x\| + \|\Im x\| = \|\Re x\|_2 + \|\Im x\|_2 \leq \sqrt{2}\|x\|_2. 
$$
Hence, if we let $\Cal H_0$ be the closure of  $\cup_n \Cal H_n \subset L^2R\subset L^2M$, 
then $p=\text{\rm proj}_{\Cal H_0}$ verifies condition $1.5$ and is thus s-smooth. 
\hfill $\square$ 

\vskip .1in 

Recall that a subset $S$  of a group $\Gamma$  is {\it free}, 
if the subgroup generated by $S$ is the free group with generators in $S$, $\Bbb F_{S}$. 

\proclaim{1.7. Proposition} With $(M,\tau)$ as in $1.4$, assume $L(\Gamma)\subset M$ is a group von Neumann subalgebra. Let $S\subset \Gamma$ 
be an infinite subset and denote by $p_S$  
the orthogonal projection of $L^2M$ onto the Hilbert space $\Cal H_S$ having orthonormal basis $\{u_g \mid g\in S\}$.

$1^\circ$ If $S$ is a free subset of $\Gamma$, 
then $p_S$ is s-smooth.

$2^\circ$ If $p_{S}$ is s-smooth, then $\Gamma$ is non-amenable. 
\endproclaim 
\noindent
{\it Proof}. $1^\circ$ By [AO] and [B], if $S\subset \Gamma$ is a free set then for any $x\in \ell^2(S)$ we have $\|x\| \leq 2\|x\|_2$. 

$2^\circ$ If $\Gamma$ is amenable, 
then by Kesten's characterization of amenability in [Ke], for any $x=\Sigma_g c_g u_g \in L(\Gamma)$ with $c_g \geq 0$ and $\Sigma_g c_g=1$, 
one has $\|x\|=1$. This shows in particular that an element of the form $x_n=n^{-1}\Sigma_{i=1}^n u_{s_i}$, with $s_i\in S$, 
has norm $\|x_n\|=1$ while $\|x_n\|_2=n^{-1/2}$. Letting $n$ tend to infinity shows that $p_S$ cannot be s-smooth. 
\hfill $\square$

\proclaim{1.8. Lemma} Let  $M$ be a von Neumann algebra,  
$\Cal B$ an operatorial Banach $M$-bimodule  and $T\in s^*(\Cal B)$ 
a smooth element.

$1^\circ$ Given any  $\varepsilon > 0$, there exists a normal state $\varphi$ on $M$ and $\alpha > 0$ such that 
if $p_1, ..., p_n \in \Cal P(M)$ is a partition of $1$ with projections that commute with $T$ and satisfy $\varphi(p_i)\leq \alpha$, $\forall i$, then $\|T\|\leq \varepsilon$. 

$2^\circ$ If $T$ commutes with a diffuse von Neumann subalgebra $B\subset M$, then $T=0$. 
\endproclaim 
\noindent
{\it Proof}. $1^\circ$ Since $T$ is smooth, given any $\varepsilon > 0$, there exists a normal state $\varphi$ on $M$ and 
$\alpha > 0$, such that if $p\in \Cal P(M)$ satisfies $\varphi(p)\leq \alpha$, then 
$\|pT\|\leq \varepsilon$. Thus, if $p_1, ..., p_n\in M$ is a partition of $1$ with projections satisfying $\varphi(p_i)\leq \alpha$, $\forall i$, 
by using the fact that $T=(\Sigma_i p_i) T = \Sigma_i p_iTp_i$, it follows that 

$$
\|T\|=\| \Sigma_i p_iTp_i \|=\max_i \|p_iTp_i\| \leq \max_i \|p_iT\|\leq \varepsilon.
$$

$2^\circ$ This is now immediate from part $1^\circ$, because $B$ diffuse implies  
that for any normal state $\varphi$ on $M$ and any $\alpha>0$ 
there exist partitions of $1$ in $B$ of $\varphi$-mesh less than $\alpha$.
\hfill $\square$

\heading 2.  Derivations  into smooth bimodules \endheading

\noindent
{\bf 2.1. Definition.} Let $\Cal B$ be a Banach $M$-bimodule. A {\it derivation} of $M$ into $\Cal B$ is a linear map $\delta: M \rightarrow \Cal B$ 
satisfying the property $\delta(xy)=x\delta(y)+\delta(x)y$, $\forall x, y \in M$. Recall from [R] that 
any derivation is automatically continuous from $M$ with the operator norm topology to $\Cal B$ with its norm topology. 
We say that $\delta$ is {\it smooth}, if it is continuous from $(M)_1$ with the $s^*$-topology to $\Cal B$ with its norm topology. 
Thus, if $M$  
is a finite von Neumann algebra with a faithful normal tracial state $\tau$, this amounts to the continuity of $\delta$ 
from $(M)_1$ with the Hilbert norm $\| \ \|_2$ given by $\tau$, to $\Cal B$ with its norm. 

If $\Cal B$ is a dual $M$-bimodule and $\delta: M \rightarrow \Cal B$ is a derivation, 
then we say that $\delta$ is {\it weakly continuous}, if it is continuous from the unit ball of $M$ with the ultraweak topology (i.e., the $\sigma(M,M_*)$-topology) to $\Cal B$ with its 
w$^*$-topology (i.e., 
the $\sigma(\Cal B, \Cal B_*)$-topology). Recall from [R]  that if the dual $M$-bimodule $\Cal B$ is normal, 
then any derivation of $M$ into $\Cal B$ is automatically weakly continuous in this sense. 

If $\delta$ is a derivation of $M$ in a Banach $M$-bimodule $\Cal B$, then we denote by $K^0_\delta$ the norm closure of 
the convex hull of $\{\delta(u)u^* \mid u\in \Cal U(M)\}$. 
If in addition $\Cal B$ is a dual $M$-bimodule, then we denote by $K_\delta$ 
the $\sigma(\Cal B, \Cal B_*)$-closure in $\Cal B$ of $K^0_\delta$. 
More generally, if $B$ is a von Neumann subalgebra of $M$, we denote by $K_{\delta, B}$ 
the $\sigma(\Cal B, \Cal B_*)$ closure in $\Cal B$ of the convex hull of $\text{\rm co} \{\delta(u)u^* \mid 
u\in \Cal U(B)\}$. Since $\|\delta(v)v^*\|\leq \|\delta\|$, $\forall u\in \Cal U(M)$, it follows that $K_\delta, K_{\delta,B}$ are $\sigma(\Cal B, \Cal B_*)$-compact subsets of 
the ball of radius $\|\delta\|$ of $\Cal B$. 

We will prove in this section that smooth valued derivations of von Neumann algebras can be ``integrated in abelian directions''.

\proclaim{2.2. Proposition}  Assume $\delta: M \rightarrow \Cal B$ is a smooth derivation of a von Neumann algebra $M$ into a Banach $M$-bimodule $\Cal B$. 
Then we have:

$1^\circ$ $\delta$ is continuous from the unit ball of $M$ with the ultraweak topology to $\Cal B$ with the $\sigma (\Cal B, \Cal B^*)$ topology 
$($i.e., $\delta$ is weakly continuous as a derivation of $M$ into the dual  $M$-bimodule $\Cal B^{**})$. 

$2^\circ$ If $\Cal B$ is a dual $M$-bimodule, then $\delta$ is weakly continuous.

$3^\circ$ If $M$ is finite then $\delta(M)\subset s^*(\Cal B)$. If in addition $\Cal B$ is a dual $M$-bimodule, 
then $K_\delta$ is a uniformly smooth subset of $s^*(\Cal B)$.  
\endproclaim
\noindent
{\it Proof}. $1^\circ$ If $\varphi: \Cal B \rightarrow \Bbb C$ is a continuous functional, 
then $(M)_1 \ni x \mapsto \varphi(\delta(x)) \in \Bbb C$ is continuous with respect to the 
$s^*$-topology, and thus also continuous with respect to the ultraweak topology on $(M)_1$.  This shows that $\delta$ is weakly continuous 
as a derivation into $\Cal B^{**}$. 

$2^\circ$ follows now from $1^\circ$. 

$3^\circ$ If $M$ is finite, then a basis of neighborhoods of $0$ in the $s^*$-topology is given 
by $\| \ \|_2$-neighborhoods with respect to traces on $M$. Thus, for any $\varepsilon>0$, there exists a normal trace 
$\tau$ on $M$ and $\alpha>0$ such that if $y\in (M)_1$, $\|y\|_2 \leq \alpha$ then $\|\delta(y)\|\leq \varepsilon/2$. 
Thus, if we take an arbitrary $x\in (M)_1$ and $y \in (M)_1$ with $\|y\|_2 \leq \alpha$, then $y\delta(x) = \delta (yx) - \delta(y)x$,  
so that $\|y\delta(x)\| \leq \|\delta(yx)\|+ \|\delta(y)x\| \leq \varepsilon$ (because $\|yx\|_2 \leq \alpha$). 
Similarly, $\|\delta(x)y\| \leq \varepsilon$. 

This shows that $\delta(M)\subset s^*(\Cal B)$. It also shows that  if  $\alpha_0>0$ is so that $x\in (M)_1$ 
$\|x\|_2 \leq \alpha_0$ implies $\|\delta(x)\|\leq \varepsilon_0/2$, 
then for any $u\in \Cal U(M)$ we have $\|x \delta(u)u^*\| \leq \|\delta(x)\| + \|\delta(xu)\|  
\leq \varepsilon_0$. Thus, if in addition $\Cal B$ is a dual $M$-bimodule, by the inferior semicontinuity of the norm on $\Cal B$ with respect to the  
$\sigma(\Cal B, \Cal B_*)$-topology, it follows that $\|x T \|\leq  \varepsilon_0$, $\forall T \in K_\delta$. Similarly, $\|T x \|\leq 
 \varepsilon_0$, $\forall T \in K_\delta$. Hence, $K_\delta \subset s^*(\Cal B)$ and $K_\delta$ is uniformly smooth. 

\hfill $\square$

\proclaim{2.3. Lemma}  Let $B$ be an amenable von Neumann algebra, $\Cal B$ a dual $B$-bimodule and  $\delta: B \rightarrow \Cal B$ a derivation. 
If either $\Cal B$ is normal over $B$, or if $B$ is abelian, then there exists $T \in K_{\delta}$ such that $\delta(b)=T b - bT, \forall b\in B$. 
\endproclaim
\noindent
{\it Proof}. Recall from [R] that $\delta$ is automatically norm continuous. 

For each $T \in \Cal B$ and $u\in \Cal U(M)$, we let $\Cal T_u(T)=u T u^* + \delta(u)u^*$.  Since $\Cal B$ is a dual $B$-bimodule, each 
affine transformation $\Cal T_u$ follows $\sigma(\Cal B, \Cal B_*)$-continuous on $\Cal B$. Noticing that  
if $v\in \Cal U(M)$, then $\Cal T_u(\delta(v)v^*)=\delta(uv)(uv)^*$, by 
taking convex combinations of elements of the form $\delta(v)v^*$, and then $\sigma(\Cal B, \Cal B_*)$-closure, this shows that $\Cal T_u(K_\delta)\subset K_\delta$ 
$\forall u\in \Cal U(B)$. Moreover, we clearly have $T_{u_1} \circ T_{u_2}=T_{u_1u_2}$, $\forall u_1, u_2 \in \Cal U(B)$. 

Since $B$ is amenable, there exists an amenable subgroup $\Cal U_0\subset \Cal U(B)$, such that 
$\Cal U_0''=B$. Moreover, if $B$ is abelian, then we can take $\Cal U_0=\Cal U(B)$. 
By the Markov-Kakutani fixed point theorem applied to the amenable group of affine transformations $\{\Cal T_u  \mid u\in \Cal U_0\}$ of the bounded 
$\sigma(\Cal B, \Cal B_*)$-compact convex 
set $K_{\delta}$, 
we get some $T \in K_{\delta}$ with the property that $\Cal T_u(T) =T$, $\forall u\in \Cal U_0$. Thus, $uT u^* + \delta(u)u^* = T$, $\forall u\in \Cal U_0$, implying that 
$\delta(u)=T u - u T$, $\forall u\in \Cal U_0$. So if $B$ is abelian and $\Cal U_0=\Cal U(B)$, then $\delta=\text{\rm ad}T$ on all $B$. 

If in turn $B$ is merely amenable but $\Cal B$ is assumed normal$/B$, then by [R] $\delta$ follows automatically weakly continuous, 
and by using the fact that $\text{\rm sp} \Cal U_0$ is weakly dense  in $B$ and that $\text{\rm ad}T$ is weakly continuous as well (by normality$/B$ of $\Cal B$),  
it follows that $\delta(u)=T b - bT$, $\forall b\in B$.  
\hfill $\square$ 

\proclaim{2.4. Corollary}  Let $\Cal B$ be a dual $M$-bimodule and $\delta: M \rightarrow \Cal B$ a smooth derivation.  
If $A\subset M$ is an abelian von Neumann subalgebra then there exists $T \in K_{\delta, A} \subset K_\delta \subset s^*(\Cal B) $ such that $\delta(a)=Ta - aT, \forall a\in A$. 
\endproclaim
\noindent
{\it Proof}. This is now immediate from Lemma 2.3 and Proposition 2.2. 
\hfill $\square$ 

\vskip .05in

We will now derive the main result of this section, showing that smooth valued derivations are inner on abelian and finite atomic von Neumann 
subalgebras.  The proof follows the idea and line of proof of the case the target bimodule is the ideal of compact operators, in [JP]: 

\proclaim{2.5. Theorem} Let $M$ be a von Neumann algebra, $\Cal B$ a dual  normal operatorial 
$M$-bimodule and $\delta: M \rightarrow \Cal B$ a derivation. Assume that $\delta$ takes values into a  
Banach sub-bimodule $\Cal B_0\subset s^*(\Cal B)$. 

$1^\circ$ If $A\subset M$ is an abelian von Neumann subalgebra, then there exists $T\in K_{\delta,A}\cap \Cal B_0$ such that 
$\delta(a)=\text{\rm ad} T(a)$, $\forall a\in A$. 

$2^\circ$ If $B\subset M$ is an atomic von Neumann subalgebra, then there exists a unique $T\in K_{\delta,B}$ such that 
$\delta(b)=\text{\rm ad} T(b)$, $\forall b\in B$, and this $T$ lies in $\Cal B_0$. 
\endproclaim 
\noindent
{\it Proof}.  $1^\circ$ By Lemma 2.3, there exists $T\in K_{\delta,A}\subset \Cal B$ such that $\delta(a)=Ta-aT$, $\forall a\in A$. 
Assume that $T \not\in \Cal B_0$, in other words $\|T\|_{\Cal B/\Cal B_0} > 0$.

Note that for any projection $p\in A$ we have $pT-Tp\in \Cal B_0$, and thus $\|T\|_{\Cal B/\Cal B_0}=\|pTp+(1-p)T(1-p)\|_{\Cal B/\Cal B_0}$. 
Note also that $\Cal B$ operatorial implies that 
$\Cal B/\Cal B_0$ is operatorial, i.e., $\|p T p + (1-p) T (1-p)\|_{\Cal B/\Cal B_0} =
\max \{\|p T p\|_{\Cal B/\Cal B_0}, \|(1-p)T (1-p)\|_{\Cal B/\Cal B_0}\}$, $\forall T \in \Cal B$, $p\in \Cal P(A)$.

Denote by $\Cal P$ the set of all projections $p\in A$ with the property that $\|p T p\|_{\Cal B/\Cal B_0}=\|T\|_{\Cal B/\Cal B_0}$.  
Assume $\Cal P$ has a minimal projection $e$. If $e$ is even minimal in $A$, i.e., $Ae=\Bbb Ce$, then for any $u\in \Cal U(A)$ 
we have: 
$$
e \delta(u)u^* e = \delta(eu)u^*e-\delta(e)e= z \delta(e) \overline{z}e  - \delta(e)e=0,
$$
where $z\in \Bbb T$ is so that $eu=ze$. This implies that $e T' e=0$ for any $T' \in K_{\delta, A}$, so in 
particular $eTe=0$, contradicting the fact that $e\in \Cal P$. 

Thus, $e$ is not minimal in $A$. But then, any non-zero projection $f\in Ae$ with $f\neq e$, satisfies $f, e-f \not\in \Cal P$, so 

$$
\| T \|_{\Cal B/\Cal B_0}=\| f T f +(e-f)T(e-f)\|_{\Cal B/\Cal B_0}
$$
$$
= \max \{ \|fTf\|_{\Cal B/\Cal B_0}, \|(e-f)T(e-f)\|_{\Cal B/\Cal B_0}\} < \|T\|_{\Cal B/\Cal B_0}, 
$$ 
a contradiction. 

This shows that $\Cal P$ cannot have minimal projections. Let $\Cal F$ be a maximal chain in $\Cal P$ and denote $f_0$ the 
infimum over all $f\in \Cal F$. Since $\Cal P$ has no minimal projections, $f_0 \not\in \Cal P$, i.e. $\|f_0Tf_0\|_{\Cal B/\Cal B_0}< \|T\|_{\Cal B/\Cal B_0}$. 
Since   
$$
\max \{\|(f-f_0)T(f-f_0)\|_{\Cal B/\Cal B_0}, \|f_0Tf_0\|_{\Cal B/\Cal B_0} \}
$$
$$ 
= \|(f-f_0)T(f-f_0) + f_0Tf_0\|_{\Cal B/\Cal B_0}=
\|fTf\|_{\Cal B/\Cal B_0} = \|T\|_{\Cal B/\Cal B_0},  
$$
it follows that $\|(f-f_0)T(f-f_0)\|_{\Cal B/\Cal B_0}=\|T\|_{\Cal B/\Cal B_0}$, $\forall f\in \Cal F$, 
or equivalently $\|f'Tf'\|_{\Cal B/\Cal B_0}=\|T\|_{\Cal B/\Cal B_0}$, $\forall f'\in \Cal F'$. 
Thus, $\Cal F'$ is a chain in $\Cal P$  with inf$\Cal F'=0$. 

Using this fact, we construct recursively a decreasing sequence of 
projections $f'_k \in \Cal F'$, such that $\|(f'_k - f'_{k+1})T(f'_k-f'_{k+1})\| > c$, $\forall k\geq 1$, 
where $c=(\|T\|_{\Cal B/\Cal B_0}-\|f_0Tf_0\|_{\Cal B/\Cal B_0})/2 > 0$. 

Assume 
$f'_1, ..., f'_n\in \Cal F'$ have been constructed to satisfy this property for all $1\leq k \leq n-1$. Since $\Cal F'$ is a 
chain decreasing to $0$ and $\|f'_nTf'_n\|\geq 2c$ (because $f'_n\in \Cal F'$), by  the normality of the bimodule structure on $\Cal B$ 
and the inferior semicontinuity with respect to the 
$\sigma(\Cal B, \Cal B_*)$-topology of the norm on $\Cal B$, it follows that there exists $f'\in \Cal F'$, $f' \leq f'_n$, 
such that $\|(f'_n - f')Tf'_n\|> \|f'_nTf'_n\|/2 \geq c$. Applying the same reasoning on the right hand side, if follows that there exists a projection $f''\in \Cal F'$, such that $f''\leq f'$ and 
$\|(f_n'-f')T(f_n'-f'')\|>c$.  Since $\|(f_n'-f'')T(f_n'-f'')\|=\|f_n'-f'\| \|(f_n'-f'')T(f_n'-f'')\|\geq \|(f_n'-f')T(f_n'-f'')\|>c$, the projection $f_{n+1}'=f''$ satisfies the required condition for $k=n$. 

We have thus constructed a sequence of mutually orthogonal projections $e_n = f'_n - f'_{n+1}\in A$, $n\geq 1$, such that $\|e_nTe_n\| > c$, $\forall n\geq 1$. 
Since $T\in K_{\delta, A}$, it follows that for each $n$, there exists $u_n\in \Cal U(A)$, such that $\|e_n \delta(u_n)u_n^* e_n\| > c$. Let now $u\in \Cal A$ 
be defined by $u=(1-\Sigma_n e_n) + \Sigma_n u_ne_n$. Then $u\in \Cal U(A)$, $ue_n=u_n, \forall n$, and we have 

$$
e_n \delta(u)u^*e_n= 
e_n\delta(e_nu)u_n^*e_n - e_n\delta(e_n)e_n =  e_n\delta(u_n)u_n^*e_n, 
$$
where we have used the fact that for any projection $p\in M$ we have $p\delta(p)p=0$ (because $p\delta(p)p=p\delta(p^2)p=2p\delta(p)p$). Thus, $\|e_n\delta(u)u^*e_n\| >c$, $\forall n$, 
contradicting the fact that $\delta(u)u^* \in s^*(\Cal B)$. 

Thus, the assumption $T\not\in \Cal B_0$ led us to a contradiction. This shows that $T\in \Cal B_0$, proving $1^\circ$.

\vskip .05in
$2^\circ$  Let $B=\oplus_{k\in J} B_k$, with each $B_k$ a finite dimensional factor, with matrix units $\{e^k_{ij}\}_{i,j}$. Denote by $\Cal V_k$ the group of unitaries 
in $B_k$ that are sums of elements of the form $\pm e^k_{ij}$. Also, for each finite set of indices $F\subset J$, denote $\Cal V_F=\oplus_{k\in F} \Cal V_k \oplus \Bbb C(1-s_F)$, 
where $s_F=\Sigma_{k\in F}s_k$, with $s_k=1_{B_k}$ the support projection of $B_k$.

We let $\Cal V=\cup_F \Cal V_F$, the union being over the finite subsets $F\subset J$. Note that $\Cal V$ is a locally finite group 
(thus amenable) with $\Cal V''=B$. By Lemma 2.3, there exists $T\in K_{\delta,B}$ such that $\delta=\text{\rm ad}T$ on $B$. We will prove that 
$T\in \Cal B_0$. 

Like in the proof of part $1^\circ$ above, assume $\|T\|_{\Cal B/\Cal B_0} >0$ and consider the set $\Cal P$ of all 
central projections $p\in \Cal P(\Cal Z(B))$ 
with the property that $\|pTp\|_{\Cal B/\Cal B_0}=\|T\|_{\Cal B/\Cal B_0}$. If $e\in \Cal P$ would be a minimal 
projection in $\Cal P$, then $s_k \leq e$ for some $k$. But $s_kTs_k\in \Cal B_0$. Indeed, this is because for each $u\in \Cal U(B)$ 
we have $s_k\delta(u)u^*s_k =s_k\delta(s_ku)u^*s_k-s_k\delta(s_k)s_k=s_k\delta(s_ku)u^*s_k$ and $us_k$ run over the unitary group $\Cal U(B)s_k=\Cal U(B_k)$, which is compact 
in the (operator) norm of $B$. Thus, $s_k K_{\delta,B}s_k = s_kK^0_{\delta,B}s_k\subset \Cal B_0$. But then $\|(e-s_k)T(e-s_k)\|_{\Cal B/\Cal B_0}= 
\|eTe\|_{\Cal B/\Cal B_0}=\|T\|_{\Cal B/\Cal B_0}$, contradicting the minimality of $e$. 

But if $\Cal P$ does not have a minimal projections, then one can proceed exactly like in the argument in $1^\circ$ above, to get a sequence of mutually 
orthogonal projections $e_n\in \Cal Z(B)$ such that $\|e_nTe_n\| > c>0$, $\forall n$.  This implies that for each $n$ there exists a unitary element $u_n\in \Cal U(B)$ 
such that $\|e_n\delta(u_n)u_n^*e_n\| > c $. Same way as in $1^\circ$, the unitary element $u=\Sigma_n u_ne_n + (1-\Sigma e_n)\in B$ 
will then satisfy $e_n\delta(u)u^*e_n=e_n\delta(u_n)u_n^*e_n$, so that on the one hand we have $\|e_n\delta(u)u^*e_n\|=\|e_n\delta(u_n)u_n^*e_n\| \geq c$, $\forall n$, 
but on the other hand, by smoothness of $\delta(u)u^*\in \Cal B_0$, we have $\lim_n \|e_n\delta(u)u^*e_n\| = 0$, a contradiction.  

Now, to prove the uniqueness of such a $T\in K_{\delta,B}$ implementing $\delta$ on $B$, notice that $s_kTs_k$ is in the norm closure of the convex hull 
of $\{s_k\delta(v)v^*s_k \mid v\in \Cal V_k\}$. But if we denote like in the proof of 2.3 by $T_v$ the transformation of $K=s_kK_{\delta,B}s_k$ given 
by $K \ni \xi \mapsto T_v(\xi)=v\xi v^* + s_k\delta(v)v^*s_k$,  then for each $\xi\in K$ the map 
$K\ni \xi \mapsto \pi(\xi)=\int_{\Cal V_k} T_v(\xi) \text{\rm d}\mu(v)$, where $\mu$ is the Haar measure 
on $\Cal V_k$, satisfies $\pi(T_{v'}(\xi))=\pi(\xi)$, for any $v'\in \Cal V_k$, because of the invariance of the Haar measure. By taking convex combinations and 
norm closure, this shows that $\pi$ has a single point range. Since the operator $T$ implements $\delta$, we also have $T_v(s_kTs_k)=s_kTs_k$, 
so $\pi(s_kTs_k)=s_kTs_k$, thus showing that any $T\in K_{\delta,B}$ that implements $\delta$ coincides with the single point $\pi(K)$ under $s_k$. 

\hfill $\square$ 

\vskip .05in
\noindent
{\bf 2.6. Remark.} Note that there are only two places in the above proof where we used the fact 
that the dual operatorial bimodule $\Cal B$ is normal: $(a)$ to deduce in part $1^\circ$ that if a net $p_i\in \Cal P(A)$ increases to 
some projection $p\in A$, then $\sup_i \|p_iT\| =\|pT\|$, $\sup_i \|Tp_i\|=\|Tp\|$, where $T\in \Cal B$ 
is an element in $K_\delta$ that implements $\delta$ on the abelian von Neumann subalgebra $A$; $(b)$ 
to deduce in part $2^\circ$ that since $\delta$ and $\text{\rm ad} T$ are weakly continuous, if they coincide on 
the $^*$-subalgebra $\text{\rm sp} \Cal V$, which is weakly dense in the amenable von Neumann algebra $B$, then they must coincide on all $B$. Related to this, one should note 
that in case $\Cal B$ is 
a dual operatorial $M$-bimodule (not necessarily normal), but $\delta$ is assumed smooth, then $K_\delta$ is a uniformly 
smooth subset of $s^*(\Cal B)$ (by 2.2.3$^\circ$), and since $T$ in $(a)$ and $(b)$ lies in $K_\delta$, it will still satisfy both properties. 
We will use this observation in the proof of Theorem 3.8.

\heading 3.  Vanishing  1-cohomology results 
\endheading

We show in this section 
that once a smooth valued derivation $\delta$  is implemented by a smooth element on a diffuse abelian von Neumann subalgebra  
satisfying some very weak regularity properties, then that element automatically implements $\delta$ on 
the entire von  Neumann algebra. 

We will also prove that smooth valued derivations are automatically smooth, i.e. if $\delta:M \rightarrow s^*(\Cal B)\subset \Cal B$, then $\delta$ 
is continuous from $(M)_1$ with the $s^*$-topology to $\Cal B$ with its norm topology. 

Note that these results can be formulated in terms of properties of the 1-cohomolgy of $M$ with coefficients in the smooth part of a bimodule $\Cal B$. 
Thus, if $M$ is a von Neumann algebra and $\Cal B$ is a Banach $M$-bimodule, then we denote by $Z^1(M, \Cal B)$ 
(respectively $Z^1_s(M,\Cal B)$) the space of derivations (resp. smooth derivations) of 
$M$ into $\Cal B$, by $B^1(M, \Cal B)$ (resp. $B^1_s(M,\Cal B)$) the subspace of inner derivations 
(resp. derivations implemented by smooth elements of $\Cal B$) and by $H^1(M, \Cal B)=Z^1(M, \Cal B)/B^1(M, \Cal B)$ 
(resp. $H^1_s(M, \Cal B)=Z^1_s(M,\Cal B)/B^1_s(M, \Cal B)$) the corresponding quotient space.

With this terminology, our results 
show that if $\Cal B$ is a dual normal operatorial $M$-bimodule, then $Z^1(M, s^*(\Cal B))=Z^1_s(M, \Cal B)$, $H^1(M, s^*(\Cal B))=H^1_s(M,\Cal B)$ 
($M$ arbitrary),  
and that if $\Cal B_0\subset s^*(\Cal B)$, then $H^1(M, \Cal B_0)=0$ whenever $M$ contains a diffuse abelian subalgebra that satisfies the weak regularity 
property. We will also prove that, under this same regularity condition for $M$, if $\Cal B$ is an arbitrary smooth operatorial Banach $M$-bimodule 
(not necessarily dual normal), 
then $H^1_s(M, \Cal B)=0$, i.e., any smooth derivation of $M$ into $\Cal B$ is inner.

The weak regularity property that we consider is  in the spirit of 
a similar concept considered for groups in (Definition 2.3 of [P5]; see also 1.2 in [IPP] for a related concept):

\vskip .05in 
\noindent
{\bf 3.1. Definition}. Let $B$ be a diffuse von Neumann subalgebra of a II$_1$ factor $M$. 
Consider the well ordered family of intermediate 
von Neumann subalgebras $B=B_0\subset B_1 \subset ... \subset B_\jmath \subset ...$ of $M$ 
constructed recursively, by transfinite 
induction, in the following way: $(a)$ for each $\jmath$, 
$B_{\jmath +1}$ is the von Neumann algebra generated by $v\in \Cal U(M)$ with $vB_\jmath v^*\cap B_\jmath$ diffuse; $(b)$ 
if $\jmath$ has no ``predecessor'', then $B_\jmath  = \overline{\cup_{n<\jmath} B_n}^w$. Notice that if $\imath_0$ is the first ordinal of 
cardinality $|M|$, then this family is constant from $\imath_0$ on.  Let $\imath$ be the first ordinal for which 
$B_{\imath + 1}=B_\imath$. The algebra $B_{\imath}$ is by definition the {\it wq-normalizer algebra} of $B$ in $M$.  
If the wq-normalizer of $B$ in $M$ is equal to $M$, then we say that $B$ is {\it wq-regular} in $M$. 
If the wq-normalizer of $B$ in $M$ is equal to $B$, we also say that $B$ is {\it  wq-malnormal} in $M$. 

\vskip .05in 

As in (2.3 of [P5]), we can alternatively characterize the wq-normalizer as follows: 

\proclaim{3.2. Lemma} Let $B\subset M$ be a diffuse von Neumann subalgebra and $N$ its wq-normalizer.  If an intermediate von Neumann 
algebra $B\subset Q \subset M$ satisfies the property that any $v\in \Cal U(M)$ with $vQv^*\cap Q$ diffuse, 
must lie in $Q$, then $Q$ contains $N$. Thus, the wq-normalizer of $B$ in $M$ is the smallest von Neumann subalgebra 
$N\subset M$ containing $B$ with the property that there exist no $v\in \Cal U(M)\setminus N$ with $vNv^*\cap N$ diffuse, or equivalently,  
the smallest wq-malnormal subalgebra of $M$ that contains $B$. 
\endproclaim
\noindent
{\it Proof}. Let $B_0 \subset ... \subset B_\imath=N$ be constructed by transfinite induction as in 3.1, with $B_\imath$ the first ordinal 
having the property 
that $B_{\imath + 1} = B_\imath$. Assume we have shown that $B_\jmath \subset Q$ for some $\jmath < \imath$. If $v\in \Cal U(M)$ 
satisfies $vB_\jmath v^*\cap B_{\jmath}$ diffuse, then $v\in B_{\jmath+1}$ by the definition of $B_{\jmath+1}$. But since $Q$ contains $B_\jmath$ we also 
have $vQv^*\cap  Q$ diffuse, so $v\in Q$. This shows that $B_{\jmath+1}\subset Q$. Also, if $\jmath \leq \imath$ 
has no predecessor and $B_n\subset Q$, $\forall n<\jmath$, then the weak closure of $\cup_{n<\jmath}B_n$ is contained in $Q$. Altogether, 
these facts imply that $B_\imath \subset Q$. 
\hfill 
$\square$ 

\proclaim{3.3. Lemma} Let $M$ be a von Neumann algebra, $\Cal B$ a dual normal Banach $M$-bimodule and $\delta: M \rightarrow \Cal B$ 
a derivation. If $\delta$ vanishes on two von Neumann subalgebras $N_1, N_2\subset M$, then it vanishes 
on the von Neumann algebra $N_1 \vee N_2$, generated by $N_1, N_2$. Thus, there exists a maximal von Neumann subalgebra 
of $M$ on which $\delta$ vanishes, which we will denote by $M_\delta$. 
\endproclaim
\noindent
{\it Proof}. Since if $\delta$ is equal to $0$ on $x, y\in M$ it is also $0$ on the product $xy$, it follows that $\delta$ vanishes 
on the algebra generated by $N_1, N_2$, which is a $^*$-algebra. By the weak continuity of $\delta$, it follows that $\delta$ also vanishes 
on $N_1\vee N_2$. 
\hfill
$\square$  

\proclaim{3.4. Lemma}  With $M, \Cal B, \delta$ as in $3.3$, assume $\delta(M)\subset s^*(\Cal B)$. 
Let $v\in \Cal U(M)$. 

$1^\circ$  If $vM_\delta v^* \cap M_\delta$ is diffuse, then $v\in M_\delta$. Thus, $M_\delta$ is wq-malnormal in $M$. 

$2^\circ$ Assume $M$ is of type $\text{\rm II}_1$ with a faithful normal trace $\tau$. For any $\varepsilon > 0$ there exists $\alpha > 0$ such that if $vA_0v^*\subset M_\delta$ for some 
$\alpha$-partition $A_0\subset M_\delta$, then $\|\delta(v)\|\leq \varepsilon$. 
Also, if $x=x^*\in M$ is so that $\{x\}' \cap M_\delta^\omega$ is diffuse, then $x\in M_\delta$. 
\endproclaim
\noindent
{\it Proof}. $1^\circ$ This is immediate from lemma 1.8.2$^\circ$, once we notice that if $x\in M_\delta $ satisfies $vxv^*\in M_\delta$, 
then $T=\delta(v)v^*$ commutes with $x$. 

$2^\circ$ The first part follows from 1.8.1$^\circ$, and it implies trivially  the last part. 

\hfill 
$\square$

\proclaim{3.5. Theorem} Let $M$ be a von Neumann algebra, $\Cal B$ a dual normal operatorial 
Banach $M$-bimodule and $\delta: M \rightarrow \Cal B$ 
a derivation that takes values in a Banach sub-bimodule $\Cal B_0\subset s^*(\Cal B)$, of the smooth part of $\Cal B$. 

$1^\circ$ $\delta$ is automatically smooth, with $K_\delta$ a uniformly smooth subset of $s^*(\Cal B)$. 

$2^\circ$ Let  $M=N_0\oplus N_1$ with $N_0$ finite atomic  
of $M$ and $N_1$ having no finite atomic direct summand. Let $B\subset M$ be a von Neumann subalgebra of the form $N_0\oplus B_1$, with $B_1\subset N_1$ 
a diffuse abelian von Neumann subalgebra. Then any  $T\in K_\delta$ that implements $\delta$ on $B$, must lie in $\Cal B_0$. 

$3^\circ$ Let $M=Q_0\oplus Q_1$ be the decomposition of $M$ with $Q_1$ of type $\text{\rm II}_1$ with atomic center 
and $Q_0$ having no $\text{\rm II}_1$ factor as direct summand. Given any abelian von Neumann subalgebra $A_1\subset Q_1$, 
if we denote $Q=Q_0\oplus A_1$, then there exists $T\in K_{\delta,Q}\cap \Cal B_0$ such that $\text{\rm ad}T=\delta$ on $Q$. 
\endproclaim
\noindent
{\it Proof}. $1^\circ$ Note first that by Theorem 2.5, $\delta$ is implemented by an element in $\Cal B_0\subset s^*(\Cal B)$ on any given abelian von Neumann 
subalgebra of $M$. Since $\text{\rm ad}T$ is smooth and takes values in $\Cal B_0$ for any $T\in \Cal B_0$, 
by subtracting such an ``inner'' derivation from $\delta$, it follows that given any abelian von Neumann subalgebra $A\subset M$, we may assume 
$\delta=0$ on $A$. 

If $M=M_0\oplus M_1$ with $M_0$ properly infinite and $M_1$ finite, then $M_0\simeq N_0\overline{\otimes} \Cal B(\ell^2\Bbb N)$. 
From the above, it follows that we may assume $\delta=0$ on a diffuse MASA $A_0\subset 1\otimes \Cal B(\ell^2\Bbb N)$. But such a MASA 
is regular in $\Cal B(\ell^2\Bbb N)$ and in fact, since its commutant contains $N_0$, it is regular in $M_0$ as well. So 
by 3.4.1$^\circ$ it follows that $\delta=0$ on all $M_0$. 

We are thus reduced to proving that $\delta$ is smooth on $M_1$, i.e., we may assume $M$ is a finite von Neumann algebra. 
By 2.5, by subtracting if necessary an inner derivation implemented by a smooth element, we may also assume $\delta$ vanishes 
on $\Cal Z(M)$. Let $\{p_i\}_i\subset \Cal P(\Cal Z(M))$ be a net of central projections increasing to $1$ 
such that $Mp_i$  has a faithful normal trace. We already know that for any sequence of mutually orthogonal projections $e_n\in Mp_i$ we have 
$\lim_n \|\delta(e_n)\|=0$ (this is because $\delta$ is implemented by 
elements in $s^*(\Cal B)$ on any abelian von Neumann subalgebra of $M$, in particular on the von Neumann algebra generated 
by $\{e_n\}_n$). By  Lemma A.1 in the Appendix, it follows that $\delta$ is smooth on each $Mp_i$. Note now that 
$\lim_i \|\delta(x(1-p_i))\|=0$ uniformly in $x\in (M)_1$. Indeed, for if there exists $c>0$ such that for all $i$, there exists $j>i$ 
and $x_j\in (M)_1$ with $\|\delta(x_j(1-p_j))\|\geq c>0$
then we can find an increasing sequence of indices $i_1 < i_2 ... $ in $M$ such that $\|\delta(x_{i_n}(p_{i_n}-p_{i_{n+1}}))\|\geq c/2$. 
But then $x=\Sigma_n x_{i_n}(p_{i_n}-p_{i_{n+1}})\in M$ and the mutually orthogonal projections $f_n=p_{i_n}-p_{i_{n+1}}$ 
satisfy $\|f_n\delta(x)\|\geq c/2$, 
$\forall n$, contradicting the smoothness of $\delta(x).$

\vskip .05in 

$2^\circ$ By Theorem 2.5, there exists $T_0\in K_\delta \cap \Cal B_0$ such that $\delta=\text{\rm ad}T_0$ on $B\subset M$. 
But then $\delta_0=\delta-\text{\rm ad}T_0$ still takes values in $\Cal B_0$, and it vanishes on $B$. By part $1^\circ$, 
$\delta$ is smooth so by Proposition 2.2, $T\in K_{\delta}$ is smooth. Thus, since $T-T_0$ is both smooth 
and it commutes with $B$, we have $T=T_0\in \Cal B_0$. 

\vskip .05in 

$3^\circ$ By 2.5.2$^\circ$, we may assume $\delta$ vanishes on the finite atomic part of $M$ as well as on a diffuse 
regular subalgebra of a copy of $\Cal B(\ell^2\Bbb N)$ that splits off the properly infinite summand of $M$. Arguing like in part 1$^\circ$, 
it follows that $\delta$ vanishes on $Q_0$ and on $A_1$. 
\hfill 
$\square$

\vskip .1in

\noindent
{\bf 3.6. Definition}. Let $M$ be a II$_1$ factor and $N\subset M$ a von Neumann subalgebra. We say that $N$ has {\it property} 
$(C')$ in $M$ (or that $N\subset M$ satisfies $(C')$) if the following holds true:  For any finite set 
$F\subset N$ and any $\varepsilon >0$, there exist $V\subset \Cal U(M)$ finite with the following 
properties: 
\vskip .05in

$(a)$ $\|x-E_{V''}(x)\|_2 \leq \varepsilon$, $\forall x\in F$; 

$(b)$ There exist diffuse abelian von Neumann subalgebras $A_v\subset \{v\}'\cap M^\omega$, $v\in V$, that mutually commute, 
i.e. $a a'=a'a$, $\forall a\in A_v, a'\in A_{v'}$, $v, v'\in V$. 
 
\vskip .05in 

If $M\subset M$ has $(C')$, we simply say that $M$ has property $(C')$. Note that this property for $M$ 
is weaker than property $(C)$ of $M$, as defined in [P2]. 

It is trivial by the definition that any abelian von Neumann subalgebra of $M$ has property $(C')$ in $M$. More generally, it is easy to see   
that any amenable von Neumann subalgebra $N\subset M$ has property $(C')$ in $M$.  
It has been shown in [P2] that if $M$ 
either has a Cartan subalgebra or has property $(\Gamma)$ of [MvN2], then $M$  has property $(C)$, and thus $(C')$ as well. 
Finally, note that property $(C')$ 
is hereditary: if $N_0\subset N\subset M$ 
and $N\subset M$ has $(C')$, then so does $N_0\subset M$.

\proclaim{3.7. Theorem} Let $M$ be a $\text{\rm II}_1$ factor,  $\Cal B$ a dual normal operatorial $M$-bimodule, 
$\Cal B_0\subset s^*(\Cal B)$ a Banach sub-bimodule  
and $\delta: M \rightarrow \Cal B_0$ a derivation. If $M$ has a diffuse von Neumann subalgebra $N$ which has the property $(C')$ in $M$ and is  
wq-regular in $M$, then there exists $T\in K_\delta \cap \Cal B_0$ such that $\delta(x)=\text{\rm ad} T(x)$, $\forall x\in M$. 
\endproclaim 
\noindent
{\it Proof}.  We first prove that given any finite set $F=F^* \subset N$ there exists $T=T(F) \in K_\delta$ such that $\delta(x)=\text{\rm ad}T(x)$, $\forall x\in F$, 
and thus for all $x$ in von Neumann algebra $M_F \subset M$ generated by $F$. 

By the property $(C')$ of $N\subset M$, 
for any given $\varepsilon>0$, there exists a finite set $V = V(\varepsilon)\subset \Cal U(M)$ with the property that $x\in_\varepsilon V''$ and such that 
for any $\alpha>0$ there exist mutually commuting finite dimensional abelian subalgebras $A_v\subset M$ with all atoms of trace at most $\alpha$, 
such that $\|v-E_{A_v'\cap M}(v)\|_2\leq \alpha$, $\forall v\in V$. Let $A=\vee_v A_v$ be the von Neumann algebra generated by $A_v, v\in V$,  
and let $S=S(V,\alpha)\in K_{\delta,A}$ be the (unique) element implementing $\delta$ on $A$. By Theorem 3.5, we have $\|\delta(v-E_{A_v'\cap M})\| 
\leq \alpha'$, with $\alpha' \rightarrow 0$ as $\alpha \rightarrow 0$. Since $\delta - \text{\rm ad}S$ is $A$-bimodular, if we denote by $\{e_k\}_k$ the 
minimal projections of $A$, it follows that $(\delta - \text{\rm ad}(S)) (E_{A'\cap M}(v))=\Sigma_k e_k(\delta(v) - [S,v])e_k$. Since $\delta(v)-[S,v] \in \Cal B_0$ 
is fixed in the process of taking $A_v$ (and thus $A$) with mesh $\alpha$ as small as we like, if we take $\alpha$ sufficiently small then by 1.8.1$^\circ$ we get 
$\|(\delta - \text{\rm ad}(S)) (E_{A_v'\cap M}(v))\|\leq \alpha'$, $\forall v\in V$. Altogether, we get the estimates 
$$
\|(\delta -\text{\rm ad}S)(v)\|\leq \|\delta(v-E_{A_v'\cap M})\| + \|(\delta - \text{\rm ad}(S)) (E_{A_v'\cap M}(v))\|\leq 2\alpha', \forall v\in V 
$$
Taking $\alpha_n$ so that $\alpha'\leq 2^{-n-1}$ and denoting $S_n=S(V,\alpha_n)$, we have thus constructed a sequence $S_n\in K_\delta$ with the 
property that $\|(\delta -\text{\rm ad}S_n)(v)\|\leq 2^{-n}$, $\forall v\in V$. Thus, if we let $S'$ be a weak limit point of $\{S_n\}_n$, then $S'\in K_{\delta}$ and 
$\delta(v)=\text{\rm ad}(S')$, $\forall v\in V$. Since $\delta$ and ad$S'$ are weakly continuous, this shows that $\delta=\text{\rm ad}S'$ on the whole von 
Neumann algebra generated by $V$. 

Take now $\varepsilon = 2^{-n}$ and $V_n\subset \Cal U(M)$ be a finite set satisfying $\|x-E_{V_n''}(x)\|_2 \leq 2^{-n}$, $\forall x\in F$,  and property $3.6 (b)$. 
Denote by $T_n$ the corresponding element $S'(V_n, 2^{-n}) \in K_\delta$ 
with the property that $\delta$ coincides with $\text{\rm ad}T_n$ 
on the von Neumann subalgebra $Q_n=V_n''$ satisfying $\|x-E_{Q_n}(x)\|_2 \leq 2^{-n}$, $\forall x\in F$. Let now $T$ be a weak limit point 
of $\{T_n\}_n$. Since $T_n\in K_\delta$ and $K_\delta$ is uniformly smooth, $\delta-\text{\rm ad}T_n$ are uniformly smooth and thus
$$
\|(\delta - \text{\rm ad}T_n)(x)\| = \|(\delta-\text{\rm ad}T_n)(x-E_{Q_n}(x))\|\rightarrow 0, \forall x\in F,
$$
implying that $\delta=\text{\rm ad}T$ on $F=F^*$. Hence, $\delta=\text{\rm ad}T$  on the whole von Neumann algebra generated by $F$. 

With this at hand, we finally take a weak limit point of the net $\{T(F)\}_F$, indexed by the finite subsets $F=F^*\subset N$, to obtain an element 
$T\in K_\delta$ that implements $\delta$ on all $N$. By Lemma  3.4, it follows that $T$ actually implements $\delta$ on the wq-normalizer of $N$ in $M$, 
thus on $M$. At the same time, if $A_0$ is a diffuse abelian von Neumann algebra of $M$ then by 3.5 there exists $S_0\in K_\delta \cap \Cal B_0$ 
implementing $\delta$ on $A_0$. Thus, $T-S_0$ is smooth and commutes with $A_0$, showing that $T=S_0\in K_\delta\cap \Cal B_0$. 

\hfill 
$\square$

\proclaim{3.8. Theorem} Let $M$ be a von Neumann algebra such that any $\text{\rm II}_1$ factor summand of $M$ 
has a diffuse, wq-normal subalgebra $N$ with the property $(C')$ in $M$. Let $\Cal B$ be a smooth operatorial Banach $M$-bimodule and   
$\delta:M \rightarrow \Cal B$ a smooth derivation. 
Then $\delta$ is implemented by an element  $T\in \Cal B$ with $\|T\|\leq \|\delta\|$. 
\endproclaim 
\noindent
{\it Proof}.  We can view $\delta$ as a 
(smooth) derivation of $M$ into the bidual $\Cal B^{**}$ of $\Cal B$, 
taking valuse in the closed sub-bimodule $\Cal B$ of $s^*(\Cal B^{**})$. Note that by 1.2.4$^\circ$ is a dual operatorial $M$-bimodule.  
By Proposition 2.2, since $\delta$ is smooth, $\delta: M \rightarrow \Cal B^{**}$ is weakly continuous. 

Let $A\subset N$ be a diffuse abelian von Neumann subalgebra. Thus, by Corollary 2.4, $\delta$ is implemented on $A$ 
by some $T_0\in  K_\delta\subset  \Cal B^{**}$, 
where $K_\delta$ is the $\sigma(\Cal B^{**}, \Cal B^*)$ closure in $\Cal B^{**}$ of $K^0_\delta\subset \Cal B$. 
By smoothness of $\delta$, $T_0$ is smooth in $B^{**}$ (i.e. it belongs to $s^*(\Cal B^{**})$), so by the proof of Theorem 2.5 
(see Remark 2.6), we actually have $T_0\in \Cal B$. 

From this point on, the argument 
in the proof of $3.7$ goes unchanged, to deduce that there exists $T\in K_\delta$ such that $\delta=\text{\rm ad}T$ on $N$, then 
by wq-normality, on all $M$. Since $\text{\rm adT}$ coincides with $\text{\rm ad}T_0$ on a diffuse subalgebra, by 1.8.2$^\circ$, we actually have $T=T_0\in \Cal B$.

\hfill 
$\square$

\vskip .1in 

\noindent
{\bf 3.9. Remark} As we mentioned at the beginning of this section, Theorems 3.7, 3.8 can be viewed as vanishing 1-cohomology results,  
showing that for II$_1$ factors 
$M$ the 1-cohomology with values in a closed submodule $\Cal B_0$ of the smooth part of 
a dual normal operatorial $M$-bimodule $\Cal B$, $H^1(M,\Cal B_0)$ (respectively the smooth 1st cohomology 
$H^1_s(M,\Cal B)$) vanishes as soon as $M$ satisfies some rather weak 
decomposability properties. In particular, this is the case if $M$ has  a Cartan subalgebra, more generally if $M$ has a 
diffuse amenable subalgebra whose quasi-normalizer generates $M$. As well as for non-prime factors $M$ (i.e., factors 
that can be decomposed as a tensor product of two II$_1$ factors) and for factors having property $(\Gamma)$ of Murray-von Neumann 
([MvN2]). 

The class of factors covered by Theorems 3.7, 3.8 contains many group factors $L(\Gamma)$, with $\Gamma$ 
infinite conjugacy class (ICC) groups. For instance,  all wreath product groups $\Gamma=H \wr G$, 
with $H$ non-trivial, are in this class. 
Indeed, if $H$ is finite, then $L(\Gamma)$ has $L(H^{(\Gamma)})$ as an amenable diffuse regular von Neumann subalgebra. 
If $|H|=\infty$, then $L(H)$ contains a diffuse abelian von Neumann subalgebra $A$, and any such algebra is wq-regular in $L(\Gamma)$. 
Another class of group factors covered by $3.7$ are $L(\Gamma_n)$ for $\Gamma_n=PSL(n, \Bbb Z)$, $n\geq 3$. 
This is because any such $\Gamma_n$ contains a chain of infinite abelian subgroups $H_1, ..., H_m$ that generate $\Gamma_n$ 
and so that $H_i$ commutes with $H_{i+1}$ for all $i$. Thus, $L(\Gamma_n)$ has diffuse abelian von Neumann subalgebra 
(e.g., $L(H_1)$) that are wq-regular in $L(\Gamma_n)$. It is interesting to note that in all these cases the $L^2$-cohomology of the group $\Gamma$ 
vanishes as well (cf. [G]). 

Note that $PSL(2,\Bbb Z)$ and the free groups $\Bbb F_t$, $2\leq t \leq \infty$,  do 
not have this property. In fact, by results in [Dy], the free group factors $L(\Bbb F_t)$ do not have property $(C)$ of [P2]. 
The same proof shows that they do not satisfy the weaker condition $(C')$ either. But despite the many in-decomposability properties 
one knows to prove for the free group factors (primeness, solidity, absence of diffuse amenable subalgebras with 
non-amenable normalizing algebra, etc), the fact that $M=L(\Bbb F_t)$ doesn't 
have any wq-regular diffuse abelian von Neumann subalgebra  is still an open problem. If one could find some smooth $L(\Bbb F_t)$-bimodule 
$\Cal B$ for which $H^1(L(\Bbb F_t), \Cal B)\neq 0$, then by 3.7, 3.8 it would follow that $L(\Bbb F_t)$ doesn't even have a diffuse wq-regular 
subalgebra with property $(C')$.

\heading 4.  The case $\Cal B = \Cal B(\Cal H)$ 
\endheading

We will prove in this section that if the target $M$-bimodule is a closed subspace $\Cal B_0$ of the smooth part of the algebra $\Cal B(\Cal H)$ of linear 
bounded operators on the Hilbert space $\Cal H$ on which 
the von Neumann algebra $M$ acts, then any derivation with values in $\Cal B_0$ is implemented by an element in $\Cal B_0$. 
In the particular case $\Cal B_0=\Cal K(\Cal H)$, this result amounts to  (Theorem II in [P3]). However, as we have seen in Section 1, there is a large 
fauna of non-compact operators in $\Cal B(\Cal H)$ that are smooth over $M$, even when $M$ is diffuse abelian. 

One should mention that, while the proof of this result 
follows closely the ideas and line of proof in [P3], the arguments in [P3] do not simply extend, as such,  to this 
larger degree of generality, and we will need to resolve additional technical difficulties, 
notably in the proof of the analogue of (Lemma 1.1 in [P3]). We overcome this through a careful usage of the 
incremental patching techniques in [P4,6,7].

\proclaim{4.1. Theorem} Let $M\subset \Cal B(\Cal H)$ be a normal representation of a von Neumann algebra 
$M$ on a Hilbert space $\Cal H$. Let $\Cal B_0\subset s^*(\Cal B(\Cal H))$ be a 
Banach $M$ sub-bimodule and $\delta: M \rightarrow \Cal B_0$ a derivation. Then there exists $T\in K_\delta \cap \Cal B_0$ such that $\delta=\text{\rm ad}T$. 
\endproclaim 
\noindent
{\it Proof}. By Theorem 3.5, we only need to prove 4.1 in the case $M$ is of type II$_1$ with atomic center and $\delta$ vanishing on $\Cal Z(M)$. 
This means, we are actually reduced to proving the case when $M$ is a II$_1$ factor. By writing $M$ as an increasing union of a net of separable II$_1$ factors 
(see e.g. the proof of A.1.2 in [P8]), 
and taking into account that if a derivation $\delta: M \rightarrow \Cal B_0\subset s^*(\Cal B(\Cal H))$ is implemented by elements $T_1, T_2\in \Cal B_0$ 
on subfactors $N_1\subset N_2 \subset M$, then $T_1, T_2$ must coincide (by $1.8.2^\circ$), it follows that it is sufficient to prove the statement in the case $M$ 
is a separable II$_1$ factor. 

By ([P1]), $M$ contains  a hyperfinite II$_1$ subfactor $R\subset M$ with trivial relative commutant, $R'\cap M=\Bbb C1$. 
By Theorem 3.7, we may assume $\delta$ vanishes on $R$. We want to prove that $\delta$ must then vanish on all $M$. We will do this 
by contradiction. If $\delta \neq 0$, then there exists $v\in \Cal U(M)$ 
such that $\delta(v)\neq 0$. Thus, there exists $\xi_0\in \Cal H$ such that the orthogonal projection $p'\in M'\cap \Cal B(\Cal H)$ of $\Cal H$ onto $\overline{M\xi_0}$ 
satisfies $p'\delta(v)p'\neq 0$. Since we also have $p's^*(\Cal B(\Cal H))p'\subset s^*(\Cal B(p'(\Cal H)))$ (by 1.3.3$^\circ$), this shows that it is sufficient to derive the 
contradiction in the case $\Cal H=L^2M$, with $M\subset \Cal B(L^2M)$ the 
standard representation of $M$. Since $s^*(\Cal B(\Cal H))$ is a $M'$-bimodule as well, by replacing if necessary $\delta$ by 
a derivation  $x\mapsto x_1'\delta(x)x_2'$ for some appropriate $x_1', x_2'\in M'$, it follows that we may actually assume $\langle \delta(v)(\hat{1}), \hat{v} \rangle =1$.  

At this point, we need some versions of (Lemma 1.1, Lemma 1.2 in [P3]), 
for $s^*(\Cal B(L^2M))$ in lieu of $\Cal K(L^2M)$. The proof of the first of these lemmas follows the line of proof of (1.1 in [P3]), using results from ([P4]) and the 
incremental patching techniques developed in ([P3, P4, P6, P7]).

\proclaim{4.2. Lemma} Assume $\Cal H$ is a separable Hilbert space and $R\subset M$ is a hyperfinite $\text{\rm II}_1$  
factor with $R'\cap M=\Bbb C$. For any countable set $\{T_m\}_m$ in the unit ball of $s_M^*(\Cal B(\Cal H))$,   
there exist unitary elements $u_n\in R$ such that $\lim_n \|T_m(u_nv)^k (\xi)\|_{\Cal H}=0$, for all $k \neq 0, \xi \in \Cal H$, $v\in \Cal U(M)$, $m\geq 1$. 
\endproclaim 
\noindent
{\it Proof}. Note first that it is sufficient to prove the statement in the case $\Cal H=L^2M$. 
Indeed, the statement only concerns spaces of the form $M(\xi)\subset \Cal H$, which are cyclic representations of $M$, thus included in the left 
regular representation. Also, if the $p_\xi$ is the orthogonal projection of $\Cal H$ onto the closure $\Cal H_\xi$ of this space, 
then we clearly have $p_\xi T_m p_\xi \in s^*(\Cal B(\Cal H_\xi))$. 

Moreover, in order to prove the statement in the case $\Cal H=L^2M$, it is clearly sufficient to 
prove that for any finite sets $F\subset L^2M$, $V\subset \Cal U(M)$, $\Cal T\subset s^*(\Cal B(L^2M))$, any $\varepsilon > 0$ and any $n\geq 1$, 
there exists a unitary element $u\in R$ such that $\| T((uv)^k\xi)\|^2 \leq \varepsilon $, $\forall \xi \in F$, 
$v\in V$, $T\in \Cal T$, $1\leq |k|\leq n$. Indeed, for if we have this, then we take $\{\xi_n\}_n$ dense in the unit ball of $L^2M$, $\{v_n\}_n$ dense the norm $\| \ \|_2$ in $\Cal U(M)$,  
and for each $n\geq 1$ we apply it to $F=\{\xi_1, ..., \xi_n\}$, $V=\{v_1, ..., v_n\}$, $\Cal T=\{T_1, ..., T_n\}$, $\varepsilon =2^{-n}$, to get a unitary element $u_n$ 
that satisfies 

$$
\| T_l((u_nv_i)^k(\xi_j)\|^2 \leq 2^{-n}, \forall i, j, |k|, l \leq n. \tag 1
$$
The sequence $\{u_n\}_n$ will then clearly satisfy the condition in the statement of the lemma, 
by the density of $\{\xi_n\}_n$ in $(L^2M)_1$ and  $\{v_n\}_n$ in $\Cal U(M)$. 
 
Denote by $\Cal W$ the set of partial isometries $w\in R$ with the properties  that $ww^*=w^*w$ and $\| T((wv)^k\xi)\|^2 \leq \varepsilon \tau(ww^*)$, $\forall \xi \in F$, 
$v\in V$, $T\in \Cal T$, $1\leq |k|\leq n$ (where $w^{-k}=(w^*)^k$ for $k>0$). 
We endow $\Cal W$ with the order given by $w_1\leq w_2$, $w_2w_1^*w_1 = w_1$. Then $(\Cal W, \leq)$ is clearly inductively ordered and we let $w$ be a maximal element. 
All we need is to prove that $w$ is a unitary element. Assume  $p=1-ww^*\neq 0$. If $w_0\in pRp$ is a partial isometry with $w_0w_0^*=w_0^*w_0=q$ 
and we denote $u=w+w_0$, then we have:
 
$$
\langle T^*T((uv)^k(\xi)), (uv)^k(\xi) \rangle = \langle T^*T((w+w_0)v)^k(\xi), ((w+w_0)v)^k(\xi) \rangle \tag 2
$$
$$
= \langle T^*T((wv)^k(\xi)), (wv)^k(\xi) \rangle + \Sigma'_{(i,0)} + \Sigma'_{(0,j)}+\Sigma'_{(i,j)}, 
$$
where the last line of $(2)$ represents the sum of $2^k \times 2^k$ elements, coming from 
developing the binomial powers  $((w+w_0)v)^k$ in the scalar product of the first line of $(2)$. 
The summations  $\Sigma'_{(i,0)}$, $\Sigma'_{(0,j)}$, $\Sigma'_{(i,j)}$ are indexed over $i, j \geq 1$ and they have the following significance: 

For each $i\geq 1$, 
$\Sigma'_{(i,0)}$ is the sum of all terms with $i$ appearances of $w_0$ on the left hand side of the scalar product and no appearance 
on the right hand side, i.e. each such term is of the form 
$$
\langle T^*T (wv)^{m_0} (\Pi_{r=1}^i w_0v(wv)^{m_r}) (\xi), (wv)^k(\xi) \rangle, \tag 3 
$$
for some $m_0, m_i \geq 0$ 
and $m_1, ..., m_{i-1}\geq 1$. 
Similarly, $\Sigma'_{(0,j)}$ is the sum of all terms with $j$ appearances of $w_0$ on the right hand side of the scalar product and no appearance 
on the left hand side, i.e. each such term is of the form 
$$
\langle T^*T (wv)^k(\xi), (wv)^{n_0}(\Pi_{s=1}^j w_0v(wv)^{n_s}) (\xi),  \rangle, \tag 4 
$$
for some $n_0, n_j \geq 0$,  $n_1, ..., n_{j-1}\geq 1$.  

For $i, j \geq 1$, $\Sigma'_{(i,j)}$ is the sum of all terms with $i$ appearances of $w_0$ on the left hand side and $j$ appearances of $w_0$ 
on the right hand side of the scalar product, i.e., each such term is of the form 
$$
\langle T^*T (wv)^{m_0} (\Pi_{r=1}^i w_0v(wv)^{m_r}) (\xi), (wv)^{n_0}(\Pi_{s=1}^j w_0v(wv)^{n_s}) (\xi)  \rangle, \tag 5 
$$
for some $m_0, m_i, n_0, n_j \geq 0$, $m_1, ..., m_{i-1}$, $n_1, ..., n_{j-1}\geq 1$. 

We will show that we can make the choice of the partial isometry $w_0\neq 0$ 
so that when estimating all the terms in $(3)$, $(4)$ and $(5)$ they add up to a quantity $\leq \varepsilon \tau(w_0^*w_0)$. 
We construct $w_0$ by first choosing its support projection $q\in R$ so that all the terms in $(5)$ are 
small, then choose a Haar unitary $w_0 \in qRq$ so that the terms in $(3), (4)$ are small. Let $\delta > 0$, which we will 
take  to be sufficiently small, depending on $\varepsilon > 0$ and $2^n$. 

\vskip .1in 

{\it Estimating $\Sigma'_{(i,j)}$}. Each one of the terms  in $(5)$ has $i\geq 1$ appearances of 
$w_0$ on the left and $j\geq 1$ on the right hand side of the scalar product. If either $i$ or $j$ is at least $2$, 
then the corresponding side of the scalar product will be of the form $ x_1 w_0 x_2  ...  w_0 x_p \xi$, 
for some finite number of possible $x_i\in (M)_1$. By using the fact that $(pRp)'\cap pMp=\Bbb Cp$ and (see Theorem A.1.2 in [P8], or Theorem 2.1 in [P4]; 
see also Theorem 0.1 in [P6]), 
it follows that there exists  
$q\in pRp$ such that for all the $x_i$ that appear this way and all $\xi'$ in a finite subset of $(pL^2Mp)_1$ we have

$$
\|qx_iq-\tau(px_ip)q\|_2\leq \delta \|q\|_2; \tag 6 
$$
$$
| \|q\xi'\|_2- \|q\|_2\|\xi'\|_2|\leq \delta \|q\|_2.
$$
\vskip .05in 
By using repeatedly the first inequality in $(6)$ combined with the triangle inequality for $\| \ \|_2$, with the Cauchy-Schwartz inequality and 
the second ineqality in $(6)$, it follows that each element in $(5)$ is at distance no more than $n^2\delta\|q\|_2^2=n^2\delta\tau(q)$ to 
an element of the form $\langle qKq (w_0^r \xi'), w_0^s \eta') \rangle$ for some $r, s \geq 1$, $K\in s^*(\Cal B(L^2M))$, $\xi', \eta' \in (pL^2Mp)_1$, 
with the number of possibilities for these terms depending on $n$.  

But the operatorial norm $\|qKq\|$ is small for $q$ having sufficiently small 
trace, by the smoothness of $K$. Also, by the second inequality in $(6)$, the ``size'' of the vectors $w_0^r\xi', w_0^s\eta'$ is $\|w_0^r \xi'\|\approx \|q\|_2\|\xi'\|$, 
$\|w_0^s \eta'\|\approx \|q\|_2\|\eta'\|$, with error controlled by $\delta\|q\|_2$. Altogether, this shows that all terms in $(6)$ have absolute value 
majorized by $n^2\delta \tau(q)$. Since there are at most $2^n \times 2^n=2^{2n}$ of them, if we take $\delta=\varepsilon n^{-2}2^{-2n}/4$, then 
the sum of the terms in $\Sigma'_{(i,j)}$, with $i\geq 1, j\geq 1$ and at least one of them $\geq 2$, is majorized by $\varepsilon/4$. 

We now estimate the summation $\Sigma'_{(1,1)}$, which consists of terms having exactly one occurrence  of $w_0$ on the left and one on the right 
of the scalar product. But these are of the form $\langle qKq w_0 q\xi', w_0q\eta'\rangle$, for which we have by the Cauchy-Scwartz inequality and 
the second part of $(6)$: 
$$
|\langle qKq w_0 q\xi', w_0q\eta'\rangle| \leq \|qKq\| \|q\xi'\|_2 \|q\eta'\|_2 \leq 3 \|qKq\| \|q\|_2^2=3\|qKq\|\tau(q) 
$$
Since $\|qKq\|$ can be made arbitrarily small for $\tau(q)$ small (by the smoothness of $K$), it follows that the $\Sigma_{(1,1)}$ can be majorized by 
$\varepsilon\tau(q)/4$ as well.

\vskip .1in 

{\it Estimating $\Sigma'_{(i,0)}, \Sigma'_{(0,j)}$}. Note that in the above estimates we only used the properties of the support 
projection $q$ of $w_0$. We will choose now $w_0\in \Cal U(qRq)$ to be a Haar unitary so that the sum over $i$ of all the summations $\Sigma'_{(i,0)}$ is majorized 
by $\varepsilon \tau(q)/4$ (similarly for $\Sigma'_{(0,j)}$). Indeed, each one of the terms in these sums 
is either of the form $\langle (\Pi_{r=1}^i w_0x_r) \xi', q (K^* \eta') \rangle$ or $\langle q(K \xi'), (\Pi_{s=1}^j w_0x_s) \eta' \rangle$, for some 
finite number of operators $K$ in the unit ball of $\Cal B(L^2M)$ and vectors $\xi', \eta'\in (L^2M)_1$. 

By the second inequality in $(6)$, 
if $i, j\geq 2$ then,  combining with the first inequality in $(6)$ and the Cauchy-Schwartz inequality, we get that each one of these terms can be 
perturbed by $2\delta \|q\|_2^2=2\delta \tau(q)$ to an element of the form $\langle w_0^i\xi', \eta'\rangle$ or $\langle \xi', w_0^j \eta' \rangle$. 
Moreover, each one of the terms in $\Sigma'_{(1,0)}, \Sigma'_{(0,1)}$ are also of this form, 
but with $i, j =1$. We can now take any Haar unitary $u_0\in qRq$ and use the fact that $\lim_n \langle u^n \xi', \eta' \rangle=0$, 
for any vectors $\xi', \eta'$, to obtain that for any finite set $F'\subset (L^2M)_1$ and any $\delta'>0$ 
there exists $n_0$ such that $|\langle u^n \xi', \eta' \rangle | \leq \delta' \tau(q)$, $\forall \xi', \eta'\in F'$, $|n|\geq n_0$. Thus, 
if we define $w_0=u_0^{n_0}$ then for sufficiently small $\delta'$ we get that $\Sigma'_{(i,0)} \leq \varepsilon\tau(q)/4$ 
and  $\Sigma'_{(0,j)} \leq \varepsilon\tau(q)/4$. 

\vskip .05in

Putting now together the two sets of estimates, one gets that $\Sigma'_{(i,j)}+\Sigma'_{(i,0)} + \Sigma'_{(0,j)} \leq \varepsilon \tau(q)=\varepsilon \tau(w_0^*w_0)$,  
which combined with $\langle T^*T((wv)^k(\xi)), (wv)^k(\xi) \rangle \leq \varepsilon \tau(w^*w)$ (due to $w$ being in $\Cal W$), shows that 
the last line in $(2)$ is further majorized by $\varepsilon \tau(u^*u)$. Thus, $u=w+w_0\in \Cal W$, with $w\leq u$, and since $w_0\neq 0$, we have $u\neq w$, 
contradicting the maximality of $w$ in $\Cal W$. Thus, $w$ is actually a unitary and the Lemma is proved. 
\hfill 
$\square$ 
\vskip 05in

The second technical lemma corresponds to (Lemma 1.2 in [P3]), but for $s^*(\Cal B)$ instead of $\Cal K(\Cal H)$. 
It shows that once $4.2$ holds true for a sequence of unitaries $u_n\in M_\delta$ then, as $n\rightarrow \infty$, 
the restriction of $\delta$ to the abelian 
von Neumann algebra generated by the unitary $u_nv$ tends to behave ``virtually'' like the derivation 
of $L^\infty(\Bbb T)$ into $\Cal B(L^2(\Bbb T))$  implemented by ad$P$, where $P$ is the projection of $L^2(T)$ onto $H^2(\Bbb T)$. 

\proclaim{4.3. Lemma} Let $\delta: M \rightarrow s^*(\Cal B(L^2M))$ be a derivation. Let $v\in \Cal U(M)$ and assume $\{u_n\}_n\subset \Cal U(M)$ is a 
sequence of unitary elements such that $\delta(u_n)=0$, $\forall n$ and $\lim_n \|T ((u_nv)^k\xi)\|=0$, $\forall \xi\in L^2M$, $k\neq 0$ and $T\in \{\delta(v), \delta(v)^*, 
\delta(v^{-1}), \delta(v^{-1})^*\}$. 
Then the sequence \{$\langle \delta((u_nv)^r)\hat{1}, (u_nv)^s \hat{1}\rangle \}_n$ tends to $\langle \delta(v)\hat{1}, \hat{v}\rangle$ if $r=s>0$ and to $0$ 
in all other cases. 
\endproclaim 
\noindent
{\it Proof}. Let $t>0$ be a positive integer. Since $\delta(u_n)=0$, we have 

$$
\delta((u_nv)^t)=\Sigma_{i=0}^{t-1} (u_nv)^iu_n\delta(v)(u_nv)^{t-i-1}; \tag 1
$$
$$
\delta((u_nv)^{-t})=\Sigma_{i=0}^{t-1}(u_nv)^{-i}\delta(v^{-1})u_n^{-1}(u_nv)^{-t+i+1}. \tag 2
$$

But by the assumptions in the hypothesis,  it follows that whenever $t-i-1\neq 0$, we have:  

$$
\|(u_nv)^iu_n \delta(v)(u_nv)^{t-i-1}\hat{1}\|=\|\delta(v)(u_nv)^{t-i-1}\hat{1}\| \rightarrow 0,  \tag 3
$$
and for any $t-i-1\geq 0$ we have
$$
\|(u_nv)^{-i}\delta(v^{-1})u_n^{-1}(u_nv)^{-t+i+1}\hat{1}\|=\|\delta(v^{-1})u_n^{-1}(u_nv)^{-t+i+1}\hat{1}\| \rightarrow 0.  \tag 4
$$

Thus, from $(1)$ and $(3)$, for any integer $t\geq 1$ we get: 
$$
\lim_n \|\delta((u_nv)^t)\hat{1} - (u_nv)^{t-1}u_n\delta(v)\hat{1}\| =0,  \tag 5
$$
while from $(2)$ and $(4)$ we get: 
$$
\lim_n \|\delta((u_nv)^{-t})\hat{1}\|=0. \tag 6
$$ 
Now, if $r<0$ and $s$ is arbitrary, then the statement follows immediately from $(6)$. If in turn $r>0$, then by $(5)$ we have for any integer $s$:  

$$
\lim_n (\langle \delta((u_nv)^r)\hat{1}, (u_nv)^s \hat{1} \rangle - \langle \delta(v)\hat{1}, u_n^{-1}(u_nv)^{s-r+1}\hat{1}\rangle) =0. \tag 7
$$

But if $s-r+1\neq 1$, then by hypothesis  $\lim_n \|\delta(v)^*u_n^{-1}(u_nv)^{s-r+1}\hat{1}\|=0$, so the second term in $(7)$ tends to $0$.  
Thus, if $s\neq r$, then $\lim_n \langle \delta((u_nv)^r)\hat{1}, (u_nv)^s \hat{1} \rangle =0$. While if $s=r >0$, then by $(7)$ again, we get 
$\lim_n \langle \delta((u_nv)^r)\hat{1}, (u_nv)^s \hat{1} \rangle = \langle \delta(v)\hat{1}, v  \hat{1}\rangle)$. 

The case $r=0$ is trivial, because then we have $(u_nv)^r=1$ and $\delta(1)=0$. 
 \hfill 
$\square$ 

\vskip .1in 
\noindent
{\it End of the proof of} $4.1$. From this point on, the rest of the proof of the theorem goes exactly the same way as the proof of (Theorem II in [P3]), by using Lemmas 4.2, 
4.3, in lieu of (Lemmas 1.1, 1.2 in [P3]), as well as (Lemma 1.3 in [P3]), which can be used unchanged. We include this last part of the argument, for completeness. 

Recall that we have reduce the proof to the case $M$ is a separable II$_1$ factor acting on $\Cal H=L^2M$ and $\delta:M\rightarrow \Cal B_0\subset s^*(\Cal B(L^2M))$ 
a derivation vanishing on a hyperfinite subfactor $R\subset M$ with $R'\cap M=\Bbb C1$. We want to prove that $\delta=0$ on all $M$. We assumed by contradiction 
that there exists $v\in \Cal U(M)$ with $\delta(v)\neq 0$ and we saw that we may assume $\langle \delta(v)\hat{1}, \hat{v} \rangle =1$. 

By 4.2, there exists a sequence of unitaries $u_n \in \Cal U(R)$ such that $\lim_n \|T(u_nv)^k\xi\|$ $=0$, for all $k\neq 0$, all $\xi\in L^2M$ and $T\in \{\delta(v), \delta(v)^*, 
\delta(v^{-1}), \delta(v^{-1})^*\}$, or $T\in \{p_m\}_m$, for some sequence of finite rank projections $p_m \rightarrow 1$. Note that 
$$
\lim_n \|p_m(u_nv)^k\xi\|=0, \forall m, \xi,
$$ 
implies that $\{(u_nv)^k\}_n$ tends weakly to $0$, for all $k\neq 0$. By (Lemma 1.3 in [P3]), this implies there exist Haar unitaries $v_n\in M$ such 
that $\lim_n \|v_n-u_nv\|$ $ =0$. Since $\langle \delta(v)\hat{1}, v \hat{1} \rangle =1$ 
and $\delta$ is norm continuous, by Lemma 4.3 this implies $\lim_n \langle \delta(v_n^r)\hat{1}, v_n^s \hat{1} \rangle $ is equal to $1$ if $r=s>0$ and to $0$ in all 
other cases. 

Take $A_n$ to be the von Neumann subalgebra of $M$ generated by the Haar unitary $v_n$ and denote by $e_n$ the orthogonal projection onto $L^2A_n$. 
Since $s_M^*(\Cal B(L^2M))\subset s_{A_n}^*(\Cal B(L^2M))$ and $s_{A_n}^*(\Cal B(L^2M))$ is an $A_n'$-bimodule with $e_n\in A_n'$, 
it follows that $e_n\delta(A_n)e_n \subset e_ns^*_{A_n}(\Cal B(L^2M))e_n=s^*_{A_n}(L^2A_n)$. Thus, 
$$
A_n \ni a \mapsto \delta_n(a)=e_n\delta(a)e_n\in s^*_{A_n}(L^2A_n),
$$ 
defines derivations with the property that $\|\delta_n\| \leq \|\delta\|$, $\forall n$. Moreover, since all $\delta_n$ are restrictions of $\delta$, 
which is smooth, the derivations $\{\delta_n\}_n$ are uniformly smooth, i.e., $\forall \varepsilon > 0$, $\exists \alpha>0$ such that 
for any given $n$, if $a\in (A_n)_1$ is so that $\|a\|\leq \alpha$, 
then $\|\delta_n(a)\|\leq \varepsilon$. 

Since $A_n\subset \Cal B(L^2A_n)$ are all spatially isomorphic to $L^\infty(\Bbb T)\subset \Cal B(L^2(\Bbb T))$, with the 
identification sending the Haar unitary generating $A_n$ to 
the operator $M_z$, acting on $f\in L^2(\Bbb T)$ by multiplication with the function $g(z)=z$, we may view all $\delta_n$ 
as derivations from $L^\infty(\Bbb T)$ into $s^*(\Cal B(L^2(\Bbb T))$, which are uniformly bounded in norm by $\|\delta\|$ and 
uniformly smooth. Moreover, by spatiality, $\lim_n \langle \delta_n(M_{z^r})\hat{1}, \hat{z^s} \rangle $ is equal to $1$ if $r=s>0$ and to $0$ in all other cases. 

We now fix a free ultrafilter $\omega$ on $\Bbb N$ and define $\Delta:L^\infty(\Bbb T)\rightarrow \Cal B(L^2(\Bbb T))$ to be the weak limit over $\omega$ of 
$\delta_n$, then it is easy to see that $\Delta$ is still a derivation, which satisfies $\|\Delta\|\leq \|\delta\|$ and is smooth (because $\delta_n$ are uniformly smooth). 
Moreover, we have 
$$
\langle \Delta(M_{z^r})\hat{1}, z^s \hat{1}\rangle = \lim_n \langle \delta(v_n^r)\hat{1}, v_n^s \hat{1} \rangle 
$$ 
is equal to $1$ if $r=s>0$ and to $0$ in all other cases. By using the derivation properties, this implies 
that $\langle \Delta(M_{z}) z^r, z^s\rangle$ is equal to $1$ if $r=0$, $s=1$ and to $0$ otherwise. But this means $\Delta(M_z)$ coincides 
with the the commutator $[P, M_z]$, where $P$ is the orthogonal projection of $L^2(\Bbb T)$ onto the subspace $H^2(\Bbb T)=\overline{\text{\rm sp}} \{z^k \mid k > 0\}$. 
Since both $\Delta$ and $\text{\rm ad}P$ are derivations and are weakly continuous ($\Delta$ is even smooth), 
and they coincide on the generator $M_z$ of $L^\infty(\Bbb T)$, it follows that $\Delta=\text{\rm ad}P$ on all $L^\infty(\Bbb T)$. 

Since $\Delta$ is smooth, this implies 
$\text{\rm ad}P$ is smooth. Moreover, since ad$P$ takes compact values on the dense $^*$-subalgebra $\text{\rm sp}\{M_{z^r} \mid r\in \Bbb Z\}$ 
and is smooth, it takes compact values on all $L^\infty(\Bbb T)$. By 2.5 or 3.7, $\text{\rm ad}P = \text{\rm ad}K$ for some $K\in \Cal K(L^2(\Bbb T))$. 
Thus $P-K$ commutes with $L^\infty(\Bbb T)$, which is maximal abelian in $\Cal B(L^2(\Bbb T))$. It follows that $P-K=M_{f}$ for some $f\in L^\infty(\Bbb T)$. 
But then 
$$
\langle (P-K)z^k, z^k\rangle = \langle M_{f} z^k, z^k\rangle =\int f(z) \text{\rm d}\nu(z)
$$
for all $k\in \Bbb Z$, with the left hand side tending to $1$ as $k\rightarrow \infty$ ant to $0$ as $k\rightarrow -\infty$, while the right hand side is constant. 
This final contradiction finishes the proof.  

\hfill 
$\square$

\proclaim{4.4. Corollary} Let $M\subset \Cal B(\Cal H)$ be a normal representation of a von Neumann algebra 
$M$ on a Hilbert space $\Cal H$. If $\delta: M \rightarrow \Cal B(\Cal H)$ is a smooth derivation, then there exists 
$T\in s_M^*(\Cal B(\Cal H))$, with $\|T\|\leq \|\delta\|$, such that $\delta=\text{\rm ad}T$. 
\endproclaim

\proclaim{4.5. Corollary}  Let $M_0$ be a $C^*$-algebra with a faithful trace $\tau$ and 
$M_0\subset \Cal B(\Cal H)$ a faithful representation of $M_0$. Let $\delta:M_0 \rightarrow \Cal B(\Cal H)$ be a derivation. Assume $\delta$ is continuous from the unit 
ball of $M_0$ with the topology given by the Hilbert norm $\|x\|_2=\tau(x^*x)^{1/2}$, $x\in M_0$, 
to $\Cal B(\Cal H)$ with the operator norm topology. Then  
there exists $T\in \Cal B(\Cal H)$ such that $\delta=\text{\rm ad}T$ and $\|T\|\leq \|\delta\|$.  More precisely, if $p_0$ denotes the projection onto the 
closure $\Cal H_0$ of  $\text{\rm sp}\{\delta(x)\xi \mid \xi \in \Cal H, x\in M_0\}$, then $[p_0, M_0]=0$, $[p_0, \delta(M_0)]=0$, 
the weak closure $M$ of $M_0p_0\subset \Cal B(\Cal H_0)$ 
is a finite von Neumann algebra, $\delta$ extends to a smooth derivation of $M$ into $ \Cal B(\Cal H_0)$ and is implemented by an 
element $T\in s^*_M(\Cal B(\Cal H_0))$, with $\|T\|\leq \|\delta\|$. 
\endproclaim
\noindent
{\it Proof}. Since $y\delta(x)\xi = (\delta(yx)-\delta(y)x)\xi$, we have $M_0\Cal H_0\subset \Cal H_0$. The projection $p_0$ onto $\Cal H_0$ 
satisfies $(1-p_0)\delta(x)=0$, $\forall x\in M_0$, and it is the smallest projection with this property. For each $T=\delta(y)p_0 \in \Cal B(\Cal H_0)$ and 
$\xi\in \Cal H_0$ we have $\|xT(\xi)\|$ small if $x\in (M_0)_1$, $\|x\|_2$ small. 

Thus, if we denote by $(\pi_\tau, \Cal H_\tau)$ the GNS representation of $M_0$ corresponding to the trace $\tau$,   
then the map $\theta: \pi(M_0) \rightarrow M_0p_0\subset \Cal B(\Cal H_0)$ defined by $\theta(\pi(x))=xp_0$ is a $^*$-algebra 
morphism which is continuous from the unit ball of $\pi(M_0)$ with the $\| \ \|_\tau$-topology to $(M_0p_0)_1$ 
with the $so$-topology.  This implies that $\theta$ extends to a $^*$-morphism, still denoted $\theta$, from 
$\pi(M_0)''$ to $M=(M_0p_0)'' \subset \Cal B(\Cal H_0)$. This means that $\ker \theta$ is given by a central projection $z_0$ in $\pi(M_0)''$, 
i.e. $M\simeq \pi_\tau(M_0)''z_0$ and the rest of the statement follows from 4.1, once we notice that the extension by continuity of $\delta$ to all $M=Mp_0$ is 
smooth (and is thus smooth valued).  
\hfill 
$\square$

\heading 5.  Further comments 
\endheading

\noindent
{\it 5.1. Generalized smooth cohomology.} Note that the proof of $2.5$ still works (and thus the conclusion in $2.5$ still holds true) if we merely 
assume that $\Cal B/\Cal B_0$ is operatorial, 
instead of the (stronger) condition $\Cal B$ operatorial. 

But an even greater degree of generality for which the arguments in the proof 
of 2.5 work exactly the same way, is the following framework, inspired by the work in [PR].  Let $M$ be a von Neumann 
algebra and $\Cal B$ be a dual normal $M$-bimodule. Assume $L\subset \Cal B_*$ is a subset of the unit ball of $\Cal B_*$ 
with the property that if $\varphi \in L$ then $\varphi (x \  \cdot \ y)\in L$ for any $x, y \in (M)_1$ (so in particular $\lambda L \subset L$ for any 
$\lambda \in \Bbb C$ with $|\lambda|\leq 1$). Moreover, we assume the set $L$ is ``separating'' for $\Cal B$, i.e., 
$\forall T\in \Cal B$, $\exists \varphi\in L$ such that $\varphi(T)\neq 0$.  

For $T\in \Cal B$, we denote $\|T\|_L=\sup \{ |\varphi(T)| \mid \varphi \in L\}$. Note that 
$\| \ \|_L$ is a norm on $\Cal B$ which is majorized by $\| \ \|$ and is complete on the unit ball of $(\Cal B, \| \ \|)$. Moreover, 
$(\Cal B, \| \ \|_L)$ is a normed $M$-bimodule and we denote 
by $s^*(\Cal B, \| \ \|_L)$ its smooth part. As in 1.2.1$^\circ$, this space is easily seen  to be 
a Banach $M$-bimodule with respect to the usual norm $\| \ \|$, and its unit ball is complete in the norm $\| \ \|_L$.

Let now $\Cal B_0\subset s^*(\Cal B, \| \ \|_L)$ be a sub-bimodule with the property that the unit ball of $(\Cal B_0, \| \ \|)$ 
is complete in the norm $\| \ \|_L$. 
In addition, we will require the norm implemented by $\| \ \|_L$ on the quotient space $\Cal B/\Cal B_0$ 
by $\|T/\Cal B_0\|_{L, ess}=\sup \{\|T-T_0\|_L \mid T_0 \in \Cal B_0\}, T\in \Cal B$,  
to be {\it operatorial}, in the sense that for any $T\in \Cal B$ and any $p\in \Cal P(M)$, we have:  
$$
\|pTp + (1-p)T(1-p)\|_{L,ess} = \max \{\|pTp\|_{L,ess}, \|(1-p)T(1-p)\|_{L,ess}\}. 
$$

\proclaim{5.1.1. Theorem} With the above assumptions and notations, 
let $\delta: M \rightarrow \Cal B_0$ be a derivation. Given any amenable von Neumann subalgebra $B\subset M$, 
$\delta$ is implemented on $B$ by some $T\in K_{\delta,B}\cap \Cal B_0$. 
Moreover, if $B$ is diffuse and wq-regular in $M$, then $T$ implements $\delta$ on all $M$. And if $\Cal B$ is a von Neumann algebra 
that contains $M$ as a von Neumann subalgebra, then $T$  implements $\delta$ on all $M$ even if $A$ is not wq-regular in $M$.  
\endproclaim 

\vskip .05in 
\noindent
{\bf 5.1.2. Exemples.} $1^\circ$ If in the above we take $L$ to be the whole unit ball of $\Cal B_*$, then $\| \ \|_L$ 
coincides with the norm $\| \ \|$. 

$2^\circ$ Let $\Cal M$ be a semifinite von Neumann factor and denote by $\Cal J\subset \Cal M$ the 
{\it compact ideal space} of $\Cal M$, i.e., the set of all elements  
$T\in \Cal M$ with the property that all spectral projections of $|T|$ corresponding to intervals $[c, \infty)$ for $c>0$ are finite 
projections in $\Cal M$. Let $Tr$ be a normal semifinite trace on $\Cal M$ and take $L\subset \Cal M_*$ to be the set of 
all normal functionals on $\Cal M$ of the form $\varphi_{x,y}(T)=Tr(yTx)$, with $x, y \in \Cal M$, 
$\|x\|\leq 1$, $Tr(x^*x)\leq 1$, $Tr(y^*y)\leq 1$.  Then the corresponding norm $\| \ \|_L$ 
satisfies all the above conditions. Moreover,  
$\Cal J$ is contained in $s^*(\Cal M, \| \ \|_L)$ and the quotient norm $\| \ \|_{L,ess}$ on $\Cal M/\Cal J$  
is operatorial. Thus, we recover in this case the framework and results obtained in [PR], for which Theorem 5.1.1 provides   
a more abstract setting, with a higher degree of generalization. 

Note that if $\Cal M = \Cal B(\Cal H)$, for some Hilbert space $\Cal H$, 
then $\Cal J = \Cal K(\Cal H)$ and $\| \ \|_L$ coincides with the operatorial norm on $\Cal B(\Cal H)$. 
But for type II$_\infty$ factors $\Cal M$, there are elements $T$ in $\Cal M$ 
for which $\| T \|_L < \| T \|$. Moreover, if $M\subset \Cal M$ is a von Neumann algebra and we view $\Cal M$ 
as a dual $M$-bimodule by left/right multiplication of the elements in $\Cal M$ by elements in $M$, then  
the norm $\| \ \|$ is not smooth 
on $\Cal J$, and is not operatorial on $\Cal M/\Cal J$. While from the above we see that $\| \ \|_L$ is both 
smooth on $\Cal J$ and operatorial on the quotient $\Cal M/\Cal J$. Thus, Theorem 5.1.1 is more general than Theorem 2.5.

\vskip .1in

\noindent
{\it 5.2. Smooth $n$-cohomology}.  We propose here a definition of higher Hochschild cohomology 
groups with smooth coefficients,  $H_s^n(M, \Cal B_0)$, $n\geq 2$, which extends the 1-cohomology in Section 3 to all $n$.  
This is done in the spirit of the classical operatorial approach to Hochschild cohomology in [KR], [JKR] (see also [SS]), 
but with additional continuity requirements.

Thus, let $\Cal B$ be a dual operatorial $M$-bimodule and $\Cal B_0\subset s^*(\Cal B)$ 
a Banach $M$ sub-bimodule. Denote by $\Cal L^n_s(M, \Cal B_0)$ the space 
of $n$-linear maps $\Phi: M^n=M \times M \times ... \times M \rightarrow \Cal B_0$, 
which are separately $s^*$-norm continuous in each variable, i.e. for each $i$ and each fixed $x_1, ..., x_{i-1}, x_{i+1}, ..., x_n$, the map 
$M \ni x \mapsto \Phi(x_1, ..., x_{i-1}, x, x_{i+1}, ..., x_n)
\in \Cal B_0$ is continuous from the unit ball of $M$ with the $s^*$-topology to $\Cal B_0$ with its norm topology. We view $\{\Cal L^n_s(M, \Cal B_0)\}_n$ as a {\it chain complex} 
with {\it boundary operators} $\partial_n: \Cal L_s^n \rightarrow \Cal L_s^{n+1}$ defined by 

$$
\partial_n (\Phi)(x_1, ..., x_{n+1})
=x_1\Phi(x_2, x_3, ..., x_{n+1})
-\Phi(x_1x_2, x_3, ...., x_{n+1})
$$
$$
+.... +(-1)^{n} \Phi(x_1, x_2, ..., x_nx_{n+1}) 
+ (-1)^{n+1} \Phi(x_1, ..., x_n)x_{n+1}. 
$$
Denote by $Z^n_s(M, \Cal B_0)$ the kernel of $\partial_{n+1}$,  by $B^n_s(M, \Cal B_0)$ the image of $\Cal L^n_s$ under $\partial_n$ and 
by  $H^n_s(M, \Cal B_0)=Z^n_s(M, \Cal B_0)/B^n_s(M, \Cal B_0)$ 
the corresponding  (vector) quotient space. 

As we mentioned in Section 3, Theorems 3.7, 3.8, 4.1, can be formulated as vanishing 1-cohomology results, in the form 
$H^1(M, \Cal B_0)=H^1_s(M, \Cal B_0)=0$. Note that smoothness is automatic 
for 1-cocycles (=derivations). Such ``automatic continuity'' phenomena for 1-cocycles are often present in cohomology theories involving operator algebras. 
But that's no longer the case for the higher cohomology. We have thus required smoothness in each variable, as part of the definition of $n$-cocycles.  

Note that in the case the target smooth bimodule $\Cal B_0$ is the ideal of compact operators on the Hilbert space 
$\Cal H$ on which $M$ acts, $\Cal B_0=\Cal K(\Cal H)$, 
smoothness in each variable is in fact automatic, once separate weak-continuity is assumed, a fact that was pointed out in (Proposition 3 of [R]). 
We recall here that statement  and include a proof,  for the reader's convenience, but also in order to emphasize the crucial way in which  the assumption  
that the target bimodule is the space of compact operators is being used. This shows the difficulty of dealing with higher cohomology 
with coefficients in arbitrary smooth  bimodules (even when restricting to cocycles that are smooth in each variable).

\proclaim{5.2.1. Lemma [Ra]} Let $(M, \tau)$ be a finite von Neumann algebra, acting normally on the Hilbert space $\Cal H$. 
If $\Phi: M^n \rightarrow \Cal K(\Cal H)$ is $n$-linear, separately norm continuous and weakly continuous in each variable, 
then $\Phi$ is separately $s^*$-norm continuous on the unit ball of $M$, in each variable. In particular, $Z^n_s(M, \Cal K(\Cal H))=Z^n_w(M, \Cal K(\Cal H))$, 
$\forall n\geq 1$.
\endproclaim 
 \noindent
{\it Proof}. By A.1, it is sufficient to prove that 
if $F: M \rightarrow \Cal K(\Cal H)$ is a bounded linear map which is weakly continuous (i.e. continuous from $(M)_1$ 
with the $\sigma(M,M_*)$-topology to $\Cal K(\Cal H)$ with the $\sigma(\Cal B(\Cal H), \Cal B(\Cal H)_*)$-toplogy), 
then for any mutually orthogonal projections $e_n \in \Cal P(M)$ we have $\lim_n \|F(e_n)\|=0$. 

Assume by contradiction that there exist 
mutually orthogonal projections $e_n$'s with $\|F(e_k)\|\geq c>0$, $\forall n$. Since $F(e_k)$ are compact operators and they tend weakly to $0$, 
it follows that we can choose recursively a fast growing sequence of integers $k_1 << k_2 << ...$ such that $F(e_{k_{n+1}})$ is 
$\varepsilon_{n+1}$-supported by a finite rank projection $p_{n+1}$ with $p_k p_{n+1} = 0$, $1\leq k \leq n$. 
If $\varepsilon_n$ are taken sufficiently small,  then by weak continuity of $F$ we have that 
$F(e)=F(\Sigma_n e_{k_n})=\Sigma_n F(e_{k_n})$ must be compact, 
while being norm close to $\Sigma_n p_n F(e_{k_n}) p_n$, which is a non-compact operator  (because under each $p_n$ there exist unit vectors on which 
$p_nF(e_{k_n})p_n$ has nom $\geq c/2$). This contradiction completes the proof.  
\hfill 
$\square$ 

\vskip .1in

\noindent
{\it 5.3. Towards a ``good'' cohomology theory for $\text{\rm II}_1$ factors}.  The work in this paper 
can be viewed as a revisitation of techniques and results in [JP] and [P2,3] in light of the recent efforts ([CS], [Th], [Pe], etc) to find a ``good''  1-cohomology theory   
for II$_1$ factors $M$. By a ``good cohomology'' we mean one   
that does not always vanish and can therefore 
detect properties of II$_1$ factors such as primeness, absence of Cartan subalgebras, infinite generation, etc. 
An ideal such theory should also be calculable and in the case of group II$_1$ factors $M=L(\Gamma)$ should reflect the cohomology theory of 
the group $\Gamma$. 

One can try 
to define such a cohomology theory on weakly dense $^*$-subalgebras of $M$, as in one of the venues proposed in [CS], or as in [Pe]. But this would require showing that 
the definition of the resulting invariant (the cohomology group, or at least its ``dimension'') does not depend on the choice of the dense subalgebra of $M$.    
Solving this type of problem has run into difficulties that seem insurmountable at the moment, similar to the problem        
of showing that Voiculescu's free entropy dimension of a set of generators of $M$ is independent on the generators. 

We have thus taken the point of view that a ``good cohomology'' for $M$ should be everywhere defined. The key for a successful such cohomology 
is the choice of the target $M$-bimodules $\Cal B$. The $M$-bimodules may be purely algebraic, with no topology on them (as in  [CS]). 
But the danger in this case is that the resulting 1-cohomology may be non-zero for any II$_1$ factor. Also, in order for a cohomology to detect 
finite/infinite generation of $M$, a derivation $\delta: M\rightarrow \Cal B$ should be uniquely determined by its values on a set of generators of $M$,  
imposing that $\delta$ satisfies some continuity properties. A topological version of the Connes-Shlyakhtenko's 
cohomology  ([CS]) has been proposed by A. Thom in [Th], where the target $M$-bimodule $\Cal B$ 
is the space  ${\text{\rm Aff}}(M\overline{\otimes}M)$, of all densely defined closable operators affiliated with $M\overline{\otimes}M$, and $\delta:M\rightarrow \Cal B$ 
are taken to be continuous from $M$ with the norm topology to $\Cal B$ with the topology of convergence in measure. 
But this 1-cohomology was shown to always vanish in [PV] (after prior work 
in this direction in [A] and [AKy]). Moreover, it was proved in [PV] that any derivation $\delta: M \rightarrow \Cal B$ that's continuous 
from $M$ with the norm topology to $\Cal B$ with any ``reasonable'' weak topology, is inner. So the approach in [CS], [Th], etc, is at a stalemate at the moment. 

But we retain from the above that the class of $M$-bimodules $\Cal B$ and derivations $\delta:M \rightarrow \Cal B$ considered should 
satisfy some topological/continuity conditions. Another important ``wishful'' feature (which in fact the cohomologies in [CS], [Th] do satisfy) is that 
derivations should be uniquely determined by their values on wq-regular diffuse von Neumann subalgebras of $M$ and, if possible, to be 
``integrable'' on abelian subalgebras. The first of this requirements is by analogy to properties in $L^2$-cohomology 
of groups (see e.g. [PeT]) and countable equivalence relations ([G]). It is a property that insures the 1-cohomology is ``small'' and calculable, and that it 
is vanishing once  
$M$ contains wq-regular abelian diffuse subalgebras.  This last property is best insured when any $\delta$ considered is inner (i.e. implemented by an element in $\Cal B$)    
on any abelian subalgebra $A\subset M$ ($\delta$ is ``integrable'' on $A$). This is not the case in the initial  Cheeger-Gromov approach to  
$L^2$-cohomology, where the $\Gamma$-module is $\Cal B=\ell^2\Gamma$, but it is the case in the Peterson-Thom recent approach in [PeT], 
where $\Cal B$ if taken to be the bigger module Aff$(L(\Gamma))$. 

The smooth operatorial bimodules and cohomology that we propose in this paper do satisfy all these requirements: $(a)$ smooth derivations 
are uniquely determined by their values on sets of generators of the von Neumann algebra; $(b)$ they are integrable on abelian subalgebras; $(c)$ they 
are uniquely determined by their values on wq-regular diffuse subalgebras. The problem is, of course, 
that any such cohomology may be vanishing. In the (optimistic!) contrary case, one needs to find a class of smooth bimodules for which 
the general results in Theorems 2.5, 3.7, 3.8 (or 0.1 in the 
introduction) hold true, but for which there do exist II$_1$ factors $M$ (e.g., $M=L(\Bbb F_n)$) that have non-inner derivations. The proof of Theorem 4.1  
gives ideas and techniques that could be used to show that all smooth cohomologies vanish. Alternatively, 4.1 provides     
``boundary conditions'' towards the search for a  ``good'' class of bimodules 
$\Cal B$, by indicating the type of conditions they should NOT satisfy: if $\Cal B$ is ``too operatorial'' 
(e.g.  $\Cal B$ a sub-bimodule in $\Cal B(\Cal H)$ like in 4.1), 
then all derivations into the smooth part of $\Cal B$ will be inner. The room to maneuver is slim, but there are ways to generalize 
even more the class of bimodules considered, while still having vanishing cohomology as in 2.5, 3.7, 3.8 (e.g., in the spirit of $5.1$).

\heading Appendix 
\endheading

We recall here, for the reader's convenience,  a result from ([P2]; see also A.1 in the Appendix of [P3]), showing that 
if  a weakly continuous linear map from a finite von Neumann algebra into a dual Banach space has the property that 
it is $s^*$-norm continuous on the the unit ball of any copy of $\ell^\infty\Bbb N$ in $M$, then it follows $s^*$-norm continuous 
on the unit ball of $M$. We include the proof from ([P2], [P3]), for the sake of completeness.

\proclaim{A.1. Lemma} Let $(M, \tau)$ be a finite von Neumann algebra 
with a faithful normal trace and $\Cal B$ a dual Banach space. Let $F:M \rightarrow \Cal B$ 
be a weakly continous linear map. Assume $\lim_n \|F(e_n)\|=0$ for any sequence of mutually orthogonal projections in $M$. 
Then $\lim_n \|F(x_n)\|=0$ for any sequence $\{x_n\}_n$ in the unit ball of $M$ that tends to $0$ in the strong operator topology. 
\endproclaim
\noindent
{\it Proof}. Let us prove first that 
$\lim_n \|F(f_n)\|=0$ for any sequence $f_n \in \Cal P(M)$ with $\tau(f_n) \rightarrow 0$. Suppose by contradiction 
that there exists such a sequence with $\|F(f_n)\|\geq c>0$, $\forall n$. By taking a subsequence, if necessary, we may assume $\Sigma_n \tau(f_n)<\infty$. 
Let $g_n = \vee_{k\geq n} f_k$ and note that $g_n$ is a decreasing sequence of projections with $\tau(g_n) \leq \Sigma_{k\geq n} \tau(f_k) \rightarrow_n 0$. 

Denote by $s_{nm}$ the support projection of $f_mg_nf_m$. Then $s_{nm}\leq f_m$ and $s_{nm}$ is majorized (in the sense of 
comparison of projections) by $g_n$. Thus, since $\tau$ is a trace, 
we have $\tau(s_{nm}) \leq \tau(g_n) \rightarrow_n 0$ for each $m$. Since $\{g_n\}_n$ is decreasing, $\{f_mg_nf_m\}_n$ is decreasing, and thus $\{s_{nm}\}_n$ 
is decreasing in $n$, for each $m$. It follows that $\{f_m - s_{nm}\}_n$ increases to $f_m$, implying that $\{F(f_m-s_{nm})\}_n$ is weakly convergent to $L(f_m)$. 
By the inferior semicontinuity of the norm on $\Cal B$ with respect to the w$^*$-topology, it follows that for each fixed $m$ and large enough $n_m$, 
we have $\|F(f_m-s_{n_m, m})\|\geq c/2$. 

This shows that we can construct recursively an increasing sequence of integers $n_1 < n_2 < ...$ such that the projections $h_k=f_{n_k}-s_{n_{k+1},n_k}$ 
satisfy $\|F(h_k)\|\geq c/2$, $\forall k$. These projections also satisfy $\tau(h_k) \leq \tau(f_{n_k}) \rightarrow_k 0$. Moreover, by the definitions of 
$h_k$ and of $s_{n_{k+1},n_k}$, it follows that $h_k g_{n_{k+1}} h_k=0$, so in particular $h_k f_{n_l}=0$ for $l\geq k+1$. Thus, $h_kh_l=0$, 
$\forall l\geq k+1$, i.e, $\{h_k\}_k$ are mutually orthogonal projections. Since $\|F(h_k)\|\geq c/2$, $\forall k$, this contradicts the fact that $F$ 
is $s^*$-norm continuous on atomic abelian von Neumann subalgebras.  

To prove that if $\{x_n\}_n$ is an arbitrary sequence in  $(M)_1$ with $\|x_n\|_2 \rightarrow 0$, then $\|F(x_n)\|\rightarrow 0$, it is clearly sufficient to show  this for $x_n=x_n^*$. 
Moreover, since $\| |x_n| \|_2 = \|x_n\|_2$, it follows that if $\|x_n\|_2 \rightarrow 0$, then $\|(x_n)_+\|_2 \rightarrow 0$, $\|(x_n)_-\|_2 \rightarrow 0$, 
showing that it is sufficient to prove the implication for sequences $x_n\in (M_+)_1$, with $\tau(x_n)\rightarrow 0$. Let $x_n=\Sigma_{m\geq 1} 2^{-n}e_{nm}$, $e_{nm}\in \Cal P(M)$, 
be the dyadic decomposition of $x_n$. It follows that $\tau(e_{nm})\rightarrow 0$, $\forall m$. Let $\varepsilon >0 $. Let $m_0$ be so that $2^{-m_0}(\|F\|+1)\leq \varepsilon/2$. 
By the first part of the proof, it follows that there exists $n_0$ such that $\|F(e_{nk})\| \leq \varepsilon/2$, $\forall n \geq n_0$ and $k\leq m_0$. 
Thus, if $n\geq n_0$ we have 

$$
\|F(x_n)\| \leq (\Sigma_{k=1}^{m_0} 2^{-k}\|F(e_{nm})\|) + (\Sigma_{k>m_0} 2^{-k}) \|F(e_{nk})\|  
$$
$$
\leq (\Sigma_{k=1}^{m_0} 2^{-k}) \varepsilon/2 + (\Sigma_{k>m_0} 2^{-k} ) \|F\| \leq \varepsilon.
$$
\hfill $\square$

\vskip .1in 
\noindent
{\bf A.2. Remark.}  If a von Neumann algebra $M$ is normally represented 
on a Hilbert space $\Cal H$ and we endow $\Cal B=\Cal B(\Cal H)$ with the $M$-bimodule 
structure given by left - right multiplication by elements in $M$, and we let $\Cal B_0\subset s^*(\Cal B)$ be the ideal of compact 
operators, $\Cal B_0=\Cal K(\Cal H)$, then the automatic smoothness of any derivation $\delta: M \rightarrow \Cal B_0$ was   
proved in [P2]. The argument consists in reducing the problem 
to the case $M\simeq \ell^\infty\Bbb N$ (via Lemma A.1), which in turn follows from the innerness of compact valued derivations on 
atomic abelian 
von Neumann algebras, established in [JP]. 

It was then noticed in [Ra] that in fact ANY weakly continuous 
linear map $F:M \rightarrow \Cal K(\Cal H)$ is $s^*$-norm continuous on atomic abelian von Neumann 
subalgebras of $M$, and thus all such weakly 
continuous linear compact valued maps are $s^*$-norm continuous on 
countable decomposable finite von Neumann algebras (cf. [P2], i.e. A.1 above). 

This general ``principle'' is of course no longer true if we replace $\Cal K(\Cal H)$ by 
the space of all smooth elements $s^*(\Cal B(\Cal H))$. Indeed, we have seen in 1.8.2$^\circ$ that if $M$ is  
a finite diffuse von Neumann algebra than the smooth part 
of $\Cal B(L^2M)$ contains infinite dimensional projections, and since $\Cal B_0=s^*(\Cal B(L^2M))$ is a hereditary 
C$^*$-algebra, any such projection $p$ satisfies $p\Cal B_0p\simeq \Cal B(L^2M)$. But there are of course plenty of weakly 
continuous linear maps from $M$ into $\Cal B(L^2M)$ which are not $s^*$-norm continuous  
on all their atomic abelian subalgebras, for instance the inclusion map $M\subset \Cal B(L^2M)$ does not satisfy 
this continuity property on any $\ell^\infty\Bbb N\hookrightarrow M$.

\head  References \endhead

\item{[AkO]} C. Akemann, Ostrand: {\it Computing norms in group $C^*$-algebras}, Amer. J. Math. {\bf 98} (1976), 1015-1047. 

\item{[A]} V. Alekseev: {\it On the first continuous $L^2$-cohomology of free group factors}, \newline  arXiv:1312.5647. 

\item{[AKy]} V. Alekseev,  D. Kyed: {\it Measure continuous derivations on von Neumann algebras and applications to $L^2$-cohomology}, arXiv:1110.6155

\item{[BGM]} U. Bader, T. Gelander, N. Monod: {\it A fixed point theorem for $L^1$ spaces},  Inventiones Math. {\bf 189} (2012), 143-148. 

\item{[Bo]} M. Bozejko: {\it On $\Lambda(p)$ sets with minimal constant in discrete noncommutative groups}, Proc. Amer. Math. Soc. 
{\bf 51} (1975), 407-412.

\item{[Ch]} E. Christensen: {\it Derivations and their extensions to perturbations of operator algebras}, in ``Proc. Sympos. Pure Math'' Vol {\bf 38}, 
Amer. Math. Soc., Providence RI,  1982, pp 261-273.

\item{[CSh]} A. Connes, D. Shlyakhtenko: {\it $L^2$-homology for von Neumann algebras}:   Reine Angew. Math. {\bf 586}  (2005), 125-168.

\item{[D]} J. Dixmier: ``Les algebres d'operateurs dans l'espace hilbertien'', Gauthier-Vill-\newline ars, Paris 1957, 1969. 

\item{[Dy]} K. Dykema: {\it Two applications of free entropy},  Math. Ann. {\bf 308} (1997), 545-558. 

\item{[G]} D. Gaboriau: {\it Invariants $\ell^2$ de r\'elations d'\'equivalence et de groupes},  Publ. Math. I.H.\'E.S. {\bf 95}
(2002), 93-150.

\item{[IPP]} A. Ioana, J. Peterson, S. Popa: {\it Amalgamated Free Products of w-Rigid Factors and Calculation of their Symmetry Groups}, 
Acta Math. {\bf 200} (2008), 85-153 (math.OA/0505589)

\item{[JKR]} B. E. Johnson, R. V. Kadison,  J. R. Ringrose: {\it Cohomology of operator algebras III: Reduction to normal cohomology}, 
Bull. Soc. Math. France, {\bf 100} (1972), 73-96.

\item{[JP]} B. Johnson, S. Parrott: {\it Operators commuting with a von Neumann algebra modulo the set of compact operators}, J. Funct. Analysis {\bf 11}  (1972),  39-61.

\item{[KR1]} R. V. Kadison,  J. R. Ringrose: {\it Cohomology of operator algebras I: Type I von Neumann algebras},  Acta Math.,  {\bf 126} (1971), 227-243. 

\item{[KR2]} R. V. Kadison,  J. R. Ringrose: {\it Cohomology of operator algebras II: Extended cobounding and the hyperfinite case},  Ark. Mat., {\bf 9} (1971), 55-63.

\item{[Ke]} H. Kesten: {\it Symmetric random walks on groups}, Trans. Amer. Math. Soc. {\bf 92} (1959), 336-354. 

\item{[Ki]} E. Kirchberg: {\it The derivation problem and the similarity problem are equivalent}, J. Operator Theory {\bf 36} (1996), 59-62. 

\item{[MvN1]} F. Murray, J. von Neumann:
{\it On rings of operators}, Ann. Math. {\bf 37} (1936), 116-229.

\item{[MvN2]} F. Murray, J. von Neumann: {\it Rings of operators
IV}, Ann. Math. {\bf 44} (1943), 716-808.

\item{[Pe]} J. Peterson: {\it $L^2$-rigidity in von Neumann algebras}, Invent. Math. {\bf 175} (2009), 417-433. 

\item{[PeT]} J. Peterson, A. Thom: {\it Group cocycles and the ring of affiliated operators}, Invent. Math. 

\item{[Pi1]} G. Pisier: ``Factorizations of linear operators and the geometry of Banach spaces'' CBMS Lecture Notes Vol. {\bf 60}, 1986. 

\item{[Pi2]} G. Pisier: ``Introduction to Operator Space Theory'', London Math. Soc. Lect. Notes {\bf 294}, 20013. 

\item{[Pi3]} G. Pisier: ``Similarity problems and completely bounded maps'', Springer Lecture Notes {\bf 1618}, 1995. 

\item{[P1]} S. Popa: {\it On a problem of R.V. Kadison on maximal
abelian *-subalgebras in factors}, Invent. Math., {\bf 65} (1981), 269-281.

\item{[P2]} S. Popa: {\it On derivations into the compacts and some
properties of type} II {\it factors}, in the ``Proceedings of XIV-th
Conference in Operator Theory'', Herculane-Timisoara 1980, I.
Gohberg (ed.), Birkhauser Verlag, 1984, pp 221-227.

\item{[P3]} S. Popa:  {\it The commutant modulo the set of compact
operators of a von Neumann algebra}, J. Funct. Analysis, {\bf
71} (1987), 393-408.

\item{[P4]} S. Popa: {\it Free independent sequences in type} II$_1$
{\it factors and related problems}, Asterisque, {\bf 232} (1995),
187-202.

\item{[P5]} S. Popa: {\it Some computations of $1$-cohomology groups
and construction of non orbit equivalent actions}, Journal of the
Inst. of Math. Jussieu {\bf 5} (2006), 309-332 (math.OA/0407199).

\item{[P6]} S. Popa: {\it A} II$_1$ {\it factor approach to the Kadison-Singer problem}, 
Comm. Math. Physics  (2014), arXiv:1303.1424

\item{[P7]} S. Popa: {\it Independence properties in subalgebras of ultraproduct}   II$_1$  {\it factors}, J. Funct. Analysis (2014), 
math.OA/1308.3982

\item{[P8]} S. Popa: {\it Classification of amenable subfactors of type} II, Acta Math. {\bf 172} (1994), 163-255. 

\item{[PRa]} S. Popa, F. Radulescu: {\it Derivations of von Neumann
factors into the compact ideal space of a semifinite algebra},
Duke Math. J., {\bf 57} (1988), 485-518.

\item{[PV]} S. Popa, S. Vaes: {\it Vanishing of the continuous first $L^2$-cohomology for} II$_1$  {\it factors}, Intern. Math. Res. Notices (2014), 
math.OA/1401.1186

\item{[PoS]} R. Powers, E. St\o rmer: {\it Free states of the canonical anticommutation relations},  
Comm. in Math. Physics {\bf 16} (1970), 1-33. 

\item{[Ra]} F. Radulescu: {\it Vanishing of} 
$H_w^2(M,\Cal K(\Cal H))$ {\it for certain finite von Neumann algebras}, Trans. Amer. Math. Soc. 
{\bf 326} (1991), 569-584.  

\item{[R]} J. Ringrose: {\it Automatic continuity of derivations of operator algebras}, J. London Math. Soc. {\bf 5} (1972), 432-438. 

\item{[SS]} A. M. Sinclair,  R. R. Smith: ``Hochschild cohomology of von Neumann algebras'', 
vol. {\bf 203},  London Math. Soc. Lecture Note Series, Cambridge University Press, Cambridge, 1995.

\item{[T]} M. Takesaki: ``Theory of operator algebras'', I, Vol. {\bf 124} of Encyclopaedia of Mathematical Sciences. 
Springer-Verlag, Berlin, 2002. Reprint of the first (1979) edition, Operator Algebras and Non-commutative Geometry, 5.

\item{[Th]} A. Thom: {\it $L^2$-cohomology for von Neumann algebras}, GAFA {\bf 18} (2008), 51-70. 

\enddocument